\documentclass[12pt,reqno]{amsart}

\usepackage[all]{xy}
\usepackage{amssymb}

\numberwithin{equation}{section}

\theoremstyle{plain}
\newtheorem{thm}{Theorem}[section]
\newtheorem*{thm*}{Theorem}
\newtheorem{cor}[thm]{Corollary}

\newtheorem{prop}[thm]{Proposition}

\theoremstyle{definition}

\def\C{\mathbb C}

\def\Z{\mathbb Z}

\def\Hom{\operatorname{Hom}}

\def\Ext{\operatorname{Ext}}
\def\Ann{\operatorname{Ann}}

\def\eps{\varepsilon}

\def\End{\operatorname{End}}
\def\Diff{\operatorname{Diff}}
\def\wh{\operatorname{Wh}}
\def\Wh{\overline{\wh}}

\def\ch{\operatorname{ch}}

\def\soc{\operatorname{soc}}
\def\socl{\overline{\soc}}
\def\rad{\operatorname{rad}}
\def\radl{\overline{\rad}}

\def\soc{\operatorname{soc}}

\def\hgt{\operatorname{ht}}

\def\a{\alpha}

\def\1{{\mathbf 1}}
\def\o{\otimes}

\def\g{\mathfrak g}
\def\h{\mathfrak h}
\def\n{\mathfrak n}

\def\ggbar{{\g \oplus \g}}

\def\<{\langle}
\def\>{\rangle}
\def\({\left(}
\def\){\right)}

\def\pd#1{{\partial_{#1}}}

\begin{document}

\title[Regular representation on the big cell and projective modules]
{Regular representation on the big cell and\\ big projective modules in the category $\mathcal O$.\\
}

\author{Konstantin Styrkas}
\address{Max-Planck-Institut f\"ur Mathematik\\ D-53111 Bonn, Germany}


\begin{abstract}

The regular representation $\mathfrak R(G_0)$ consists of the regular functions on the
big cell $G_0 \subset G$, arising from the Gauss decomposition of a simple complex Lie group $G$.
The algebra $\mathfrak R(G_0)$ is isomorphic to the subalgebra of $\mathcal U(\g)^*$, spanned
by the matrix elements of integral weight modules in the category $\mathcal O$. 

We study the blocks $\mathfrak R(G_0)_\lambda$, corresponding to 
central characters, associated with anti-dominant integral weights $\lambda$.
We show that $\mathfrak R(G_0)_\lambda$ is spanned by the matrix elements of the big
projective module $P_\lambda$, and prove that $\mathfrak R(G_0)_\lambda \cong P_\lambda^* \o_{Z(\g)} P_\lambda$,
which gives a version of the Peter-Weyl theorem for the projective modules.

The Whittaker functor associates to any $\g$-module the subspace of Whittaker vectors, transforming
under $\n_+$ according to a one-dimensional representation $\boldsymbol\eta$. We apply the Whittaker functor to the $\ggbar$-module
$\mathfrak R(G_0)_\lambda$, and prove that for nonsingular $\boldsymbol\eta$ we have $\Wh_\eta (\mathfrak R(G_0)_\lambda) \cong P_\lambda$,
which yields the big projective module analogue of the Borel-Weil realization of the simple finite-dimensional $\g$-modules.

Our results admit an invariant reformulation in terms of the algebras of endomorphisms of 
projective generators of $\mathcal O_\lambda$. This allows for immediate generalizations to the 
quantum group case, which also hold for the small quantum groups when $q$ is a root of unity.
\end{abstract}


\maketitle

\setcounter{section}{-1}
\section{Introduction.}

A fundamental object in the representation theory of a simple complex Lie group $G$ is the regular representation $\mathfrak R(G)$, 
realized as the algebra of regular functions on $G$. The space $\mathfrak R(G)$ carries compatible structures of a commutative associative algebra
and of a $G \times G$-module, which encode important information about the category of finite-dimensional representations of $G$.

The algebro-geometric version of the classical Peter-Weyl theorem asserts that $\mathfrak R(G)$ decomposes into a direct sum of subspaces 
$\mathbb M(L_\lambda)$, spanned by matrix elements of finite-dimensional $G$-modules $L_\lambda$, indexed by the set $\mathbf P^+$
of dominant integral weights. In other words, we have a decomposition
\begin{equation}
\mathfrak R(G) = \bigoplus_{\lambda \in \mathbf P^+} \mathbb M(L_\lambda) \cong \bigoplus_{\lambda \in \mathbf P^+} L_\lambda^* \o L_\lambda,
\end{equation}
where $L_\lambda^*$ are the modules dual to $L_\lambda$, and the second isomorphism of $G \times G$-modules is a reformulation of the fact that
the matrix elements, corresponding to any fixed basis of $L_\lambda$, are linearly independent and form a basis of $\mathbb M(L_\lambda)$. 

Another classical result provides realizations of simple $G$-modules inside the regular representation $\mathfrak R(G)$.
Let $B_-$ be a Borel subgroup of $G$, and let $T$ be its maximal torus. For any $\lambda \in \mathbf P^+$ we consider the
character $\xi_\lambda: B_- \to \C$, determined by $\lambda$ on $T$ and trivially extended to the unipotent part $N_-$ of $B_-$. 
The Borel-Weil theorem \cite{Bo} asserts that the simple module $L_\lambda$ can be realized as the subspace of $\mathfrak R(G)$,
\begin{equation}
\label{eq:Borel-Weil realization}
L_\lambda\  \cong \ \left\{\varphi \in \mathfrak R(G) \ \biggr| \ \varphi (b\, x) =
 \xi_\lambda(b) \, \varphi(x) \text{ for any } b \in B_-,\, x \in G\right\},
\end{equation}
on which $G$ acts by the right shifts, $(g\,\varphi)(x) = \varphi(x\, g)$ for any $x,g \in G$.
It also follows from the Borel-Weil theorem that the subspace of $N_-$-invariant elements in $\mathfrak R(G)$ 
gives a model for finite-dimensional $G$-modules, or equivalently
\begin{equation}
\label{eq:model group representations}
\mathfrak R(N_-\backslash G) \cong \left\{\varphi \in \mathfrak R(G) \ \biggr| \ \varphi (n\, g) = \varphi(g) \text{ for any } n \in N_-,\, g \in G\right\}
\cong \bigoplus_{\lambda \in \mathbf P^+} L_\lambda.
\end{equation}

Every $G$-module has a natural infinitesimal action of the Lie algebra $\g$ of the group $G$.
Elements of its universal enveloping algebra $\mathcal U(\g)$ can be identified with the differential operators acting
on germs of functions at the unit element $e \in G$, and evaluation at $e$ yields a linear map $\vartheta: \mathfrak R(G) \to \mathcal U(\g)^*$.
The image of $\mathfrak R(G)$ is the so-called Hopf dual subalgebra of $\mathcal U(\g)^*$, and the Lie algebraic version of
the Peter-Weyl theorem identifies it with the subspace spanned by the matrix elements of finite-dimensional $\g$-modules $L_\lambda$ for
all $\lambda \in \mathbf P^+$.

In this paper we consider the extension $\mathfrak R(G_0)$ of the regular representation $\mathfrak R(G)$, which consists of regular
functions on the big cell $G_0 \subset G$ of the Gauss decomposition of the group $G$. We call the algebra $\mathfrak R(G_0)$ the
regular representation on the big cell.
The group $G$ does not act on $\mathfrak R(G_0)$ because $G_0$ is not invariant under left and right shifts, but
the infinitesimal regular actions of $\g$ are still well-defined. The $\ggbar$-module structure of $\mathfrak R(G_0)$ cannot be
described in terms of finite-dimensional $\g$-modules, so we are naturally led to a more general class of representations of $\g$.

The Gauss decomposition $G_0 = N_- \, T \, N_+$ of the group $G$ induces the triangular decomposition 
$\g = \n_- \oplus \h \oplus \n_+$ of the Lie algebra $\g$.
The Bernstein-Gelfand-Gelfand category $\mathcal O$ consists of finitely generated, $\h$-diagonalizable and locally $\n_+$-finite $\g$-modules,
and contains simple highest weight modules $L_\lambda$ associated with $\lambda$ not necessarily dominant. The category $\mathcal O$ is
not semisimple, and each simple modules $L_\lambda$ has the indecomposable projective cover $P_\lambda$. 
The projective modules $P_\lambda$ corresponding to 
strictly anti-dominant $\lambda$ play a particularly important role; we call them the big projective modules.

Similarly to $\mathfrak R(G)$, the algebra $\mathfrak R(G_0)$ can be identified with a subalgebra of $\mathcal U(\g)^*$.
We prove that this subalgebra is spanned by the matrix elements of all modules with integral weights in $\mathcal O$, and decomposes into
a direct sum of $\ggbar$-submodules $\mathbb M_\lambda \subset \mathcal U(\g)^*$, 
\begin{equation}
\mathfrak R(G_0) \cong \bigoplus_{\lambda \in - \mathbf P^{++}} \mathbb M_\lambda, \qquad 
\mathbb M_\lambda = \sum_{V \in \mathcal O_\lambda} \mathbb M(V)
\end{equation}
where the subspaces $\mathbb M_\lambda$ are spanned by matrix elements 
of the modules in the block $\mathcal O_\lambda$, corresponding to the strictly anti-dominant integral weights $\lambda \in -\mathbf P^{++}$.
Furthermore, we prove that it suffices to consider the matrix elements of the big projective modules, i.e. for every $\lambda \in -\mathbf P^{++}$
we have $\mathbb M_\lambda = \mathbb M(P_\lambda)$.

In contrast to the finite-dimensional $G$-modules, the matrix elements, corresponding to a fixed basis of the big projective module $P_\lambda$,
are not linearly independent. Therefore, the $\ggbar$-module $\mathbb M(P_\lambda)$ is not isomorphic to $P_\lambda^* \o P_\lambda$, but
rather to its quotient by the ideal, corresponding to identically vanishing matrix elements. We establish a $\ggbar$-module isomorphism
$\mathbb M(P_\lambda) \cong P_\lambda^* \o_{Z(\g)} P_\lambda$, where $Z(\g)$ is the center of the universal enveloping algebra $\mathcal U(\g)$.
Unlike the finite-dimensional case, the action of $Z(\g)$ on $P_\lambda$ does not reduce to scalars, and generates the entire endomorphism ring 
$\End_{\mathcal O}(P_\lambda)$ of the big projective module \cite{So}. Thus, we obtain $\ggbar$-module isomorphisms
\begin{equation}
\mathfrak R(G_0) \cong \bigoplus_{\lambda \in -\mathbf P^{++}} \mathbb M(P_\lambda) \cong 
\bigoplus_{\lambda \in -\mathbf P^{++}} P_\lambda^* \o_{Z(\g)} P_\lambda,
\end{equation}
which gives an analogue of the Peter-Weyl theorem for $\mathfrak R(G_0)$ and the big projective modules. Another useful
description of $\mathfrak R(G_0)$ is given by
\begin{equation}
\label{eq:Koszul description}
\mathfrak R(G_0) \cong \bigoplus_{\lambda \in -\mathbf P^{++}} \mathcal P_\lambda^* \o_{\mathcal A_\lambda} \mathcal P_\lambda,
\end{equation}
where $\mathcal A_\lambda$ is the algebra of endomorphisms of the projective generator $\mathcal P_\lambda$ of the
block $\mathcal O_\lambda$.

The big projective modules $P_\lambda$ for $\lambda \in -\mathbf P^{++}$ admit a realization, similar to the one
provided by the Borel-Weil theorem for finite-dimensional $G$-modules. However, instead of the $N_-$-invariant elements of $\mathfrak R(G)$, 
we need to consider the more general notion of Whittaker vectors, studied by Kostant \cite{K}.
For any character $\boldsymbol \eta^-: N_- \to \C$, we consider
\begin{equation}
\label{eq:model big projective}
(\Wh^-_\eta  \o 1) (\mathfrak R(G_0)) = \left\{\varphi \in \overline{\mathfrak R(G_0)} \ \biggr| \ \varphi (n\, x) = \boldsymbol \eta^-(n)\, \varphi(x) \text{ for any } n \in N_-,\, x \in G_0\right\},
\end{equation}
where $\overline{\mathfrak R(G_0)}$ is the power series completion of the polynomial algebra $\mathfrak R(G_0)$.
The notation $(\Wh^-_\eta  \o 1)$ indicates that we consider the subspace of functions with special transformation properties 
under the left shifts by $N_-$, which remains a $\g$-module with respect to the right regular action.
We prove that for nonsingular $\boldsymbol\eta^-$ the Whittaker vectors in $\mathfrak R(G_0)$ give a model for big projective $\g$ modules,
corresponding to strictly anti-dominant weights:
\begin{equation}
(\Wh^-_\eta  \o 1) (\mathfrak R(G_0)) \cong \bigoplus_{\lambda \in -\mathbf P^{++}} P_\lambda,
\end{equation}
Imposing further transformation properties with respect to the center $Z(\g)$, one extracts the individual big projective modules $P_\lambda$
for each $\lambda \in -\mathbf P^{++}$.

The universal enveloping algebra $\mathcal U(\g)$ can be deformed into the quantum group $\mathcal U_q(\g)$,
and the Lie algebraic realization of $\mathfrak R(G_0)$ as a subspace of $\mathcal U(\g)^*$ yields 
the ``quantum coordinate algebra'' $\mathfrak R_q(G_0)$, defined as the appropriate subspace in
$\mathcal U_q(\g)^*$.
When $q$ is generic, the representation theories of $\mathcal U(\g)$ and $\mathcal U_q(\g)$ are completely parallel,
and the structural results obtained for $\mathfrak R(G_0)$ also hold for $\mathfrak R_q(G_0)$.
In particular, we get the quantum versions of the Peter-Weyl and Borel-Weil theorems for quantized big projective modules.

When $q$ is a root of unity, one can define a finite-dimensional version $\mathfrak U$ of the quantum group $\mathcal U_q(\g)$.
The dual algebra  $\mathfrak U^*$ exhibits properties similar to those of $\mathfrak R(G_0)$, and its structure
can be described as in \eqref{eq:Koszul description} in terms of projective generators in the appropriate category of representations, and their
endomorphism algebras.

I am deeply grateful to my teacher Igor B. Frenkel for numerous discussions and valuable advice. 
I also thank the Max-Planck-Institut f\"ur Mathematik in Bonn for their hospitality.

\section{Notation and preliminaries.}

\subsection{Representations of $\g$ and category $\mathcal O$}

Let $\g$ be a simple complex Lie algebra with a Cartan subalgebra $\h$.
We denote by $\Delta$ the root system of $\g$, and make a choice of
positive roots $\Delta_+$ and simple roots $\Pi$.

Denote $r = \dim \h$, and let $\{\mathbf e_i, \mathbf h_i, \mathbf f_i\}_{i=1}^{r}$ be
the Chevalley generators of the Lie algebra $\g$. The nilpotent subalgebras $\n_-$ and $\n_+$ are generated by 
$\{\mathbf f_i\}$ and $\{\mathbf e_i\}$ respectively. We choose root vectors
$\{\mathbf e_\beta,\mathbf f_\beta \}_{\beta\in \Delta_+}$ 
such that for any simple root $\a_i$ we have 
$\mathbf e_{\alpha_i} \equiv \mathbf e_i$ and $\mathbf f_{\alpha_i} \equiv \mathbf f_i$.

We refer to elements of $\h^*$ as weights. A weight $\lambda$ is called integral (resp. dominant and strictly dominant)
if for all $i=1,\dots,r$ we have $\<\lambda,\mathbf h_i\> \in \Z$  (resp. $\<\lambda,\mathbf h_i\> \ge 0$ and $\<\lambda,\mathbf h_i\> > 0$).
The sets of integral (resp. dominant integral and strictly dominant integral) weights are denoted $\mathbf P$ (resp. $\mathbf P^+$ and $\mathbf P^{++}$).
Elements of $-\mathbf P^+$ (resp. $-\mathbf P^{++}$) are called antidominant (resp. strictly antidominant) weights.

For any Weyl group element $w \in W$, the length $l(w)$ is defined as the smallest number $l$ such that
$w$ can be represented as $w = s_{j_1}\dots s_{j_l}$, where $\{s_i\}_{i=1}^r$ are the simple root reflections.
We denote $w_\circ$ the unique longest element of $W$.

We set $\rho = \frac 12 \sum_{\beta\in\Delta_+} \beta$, and define the ``dot'' action of the Weyl group on $\h^*$ by
\begin{equation}
s_i\cdot\lambda = \lambda - \<\lambda+\rho,\mathbf h_i\> \a_i.
\end{equation}
For any $\lambda \in \h^*$, denote $W_\lambda$ the subgroup of $W$ stabilizing $\lambda$, and
pick a set $W^\lambda$ of representatives of the cosets $W/W_\lambda$. We call $\lambda$ regular if $W_\lambda = \{e\}$.
We also denote by $l_\lambda$ the smallest length of an element $w \in W$ such that $w\cdot\lambda = w_\circ\cdot\lambda$; 
if $\lambda$ is regular, then $l_\lambda = l(w_\circ) = |\Delta_+|$.

The Bernstein-Gelfand-Gelfand category $\mathcal O$ consists of finitely generated, locally $\n_+$-finite
$\g$-modules $V$, decomposing into the direct sum of finite-dimensional weight subspaces
\begin{equation}
V = \bigoplus_{\mu\in\h^*} V[\mu], \qquad V[\mu] = 
\{ v \in V | \mathbf h \,v = \<\mu,\mathbf h\> \, v \text{ for any } \mathbf h \in \h\}.
\end{equation}
The formal character of $V$ is defined by
\begin{equation}
\label{eq:formal character}
\ch V = \sum_\mu \dim V[\mu] \ e^\mu.
\end{equation}

Simple objects in $\mathcal O$ are the irreducible highest weight modules $L_\lambda$, parameterized by $\lambda \in \h^*$.
Each simple module $L_\lambda$ has a projective cover $P_\lambda$ with simple top, isomorphic to $L_\lambda$. The category
$\mathcal O$ also contains the ``standard'' Verma module $M_\lambda$ and the ``costandard'' contragredient Verma module $M_\lambda^c$
for each $\lambda\in \h^*$.

A central character of $\g$ is a homomorphism $\chi:Z(\g) \to \C$, where $Z(\g)$ is the center of the 
universal enveloping algebra $\mathcal U(\g)$. For any central character $\chi$, we define 
\begin{equation}
\label{eq:category O blocks}
\mathcal O(\chi) = \left\{ V \in \mathcal O \, \bigr| \, \forall z  \in Z(\g) \text { the operator } z-\chi(z) \text{ acts nilpotently in } V \right \}
\end{equation}
Each $\mathcal O(\chi)$ is a full subcategory of $\mathcal O$, and we have the decomposition $\mathcal O = \bigoplus_\chi \mathcal O(\chi)$, such that $\Ext_{\mathcal O}^\bullet(M,N) = 0$ if $M \in \mathcal O(\chi), N\in\mathcal O(\chi')$ and $\chi \ne \chi'$.

For any $\lambda \in \h^*$ elements $z \in Z(\g)$ act in the simple module $L_\lambda$ as scalars. We denote $\chi_\lambda$ the corresponding central character; then $\chi_\lambda = \chi_\mu$ if and only if $\lambda = w\cdot\mu$ for some $w \in W$. 

We denote $\mathcal O_\lambda = \mathcal O(\chi_\lambda)$.
For generic $\lambda \in \h^*$ the simple module $L_\lambda$ is projective, and the category $\mathcal O_\lambda$ consists of
finite direct sums of $L_\lambda$. In this paper we study the blocks $\mathcal O_\lambda$ associated with integral weights $\lambda$,
which have much richer structure.

All projective modules in $\mathcal O_\lambda$ have Verma flags, i.e. filtrations with quotients are isomorphic to Verma modules.
Denote $(P_{x\cdot\lambda}:M_{y\cdot\lambda})$ the number of times $M_{y\cdot\lambda}$ appears as a quotient in a Verma flag
of $P_{x\cdot\lambda}$; it does not depend on the choice of the filtration.
The Bernstein-Gelfand-Gelfand reciprocity states that for any $x,y \in W^\lambda$
\begin{equation}
(P_{x\cdot\lambda}:M_{y\cdot\lambda}) = [M_{y\cdot\lambda}:L_{x\cdot\lambda}],
\end{equation}
where $[M_{y\cdot\lambda}:L_{x\cdot\lambda}]$ denotes the multiplicity of the simple modules $L_{x\cdot\lambda}$ in
the composition series for $M_{y\cdot\lambda}$, computed by the appropriate Kazhdan-Lusztig polynomial.
As a consequence, one obtains the symmetry of the so-called Cartan matrix of $\mathcal O_\lambda$:
\begin{equation*}
[P_{x\cdot\lambda}:L_{y\cdot\lambda}] = \sum_{w\in W^\lambda} [P_{x\cdot\lambda}:M_{w\cdot\lambda}][M_{w\cdot\lambda}:L_{y\cdot\lambda}] = 
\sum_{w\in W^\lambda} [M_{w\cdot\lambda}:L_{x\cdot\lambda}][P_{y\cdot\lambda}:M_{w\cdot\lambda}] = [P_{y\cdot\lambda}:L_{x\cdot\lambda}].
\end{equation*}

Let $\lambda$ be a dominant integral weight. There are two important modules in the block $\mathcal O_\lambda$: the unique
finite-dimensional simple module $L_\lambda$, and the big projective module $P_{w_\circ\cdot\lambda}$. One of the
remarkable similarities between them is exhibited by via their characters.
The Weyl character formula for the finite-dimensional module $L_\lambda$ can be written as
\begin{equation}
\label{eq:Weyl character}
\ch L_\lambda = \sum_{w\in W^\lambda} (-1)^{|w|} \ch M_{w\cdot\lambda},
\end{equation}
and we have a similar expression for the character of the big projective module:
\begin{equation}
\label{eq:big projective character}
\ch P_{w_\circ\cdot\lambda} = \sum_{w\in W^\lambda} \ch M_{w\cdot\lambda}.
\end{equation}

For any $V \in \mathcal O$ we define the restricted dual $V^* = \bigoplus_\mu V[\mu]^*$.
Dualizing the $\g$-action on $V$, we get a natural {\it right} action of $\g$ on $V^*$,
such that $\<v^*, x\, v\> = \< v^*\, x,  v\>$ for any $x \in \g$ and $v \in V, v^* \in V^*$.
We transform it into the more traditional left $\g$-action by means of the Lie algebra
anti-involution $x \mapsto x^* = -x$.
In other words, we have
\begin{equation}
\label{eq:antipode involution}
\<x \, v^*,  v\> = \< v^* x^*,  v\> = - \< v^* \, x,  v\> = - \< v^*, x \, v\>, \qquad x \in \g.
\end{equation}
The dual module $V^*$ belongs to the ``mirror'' category $\mathcal O^*$, associated with the lowest weight $\g$-modules. 
In other words, $\mathcal O^*$ consists of finitely generated, $\h$-diagonalizable and locally $\n_-$-nilpotent $\g$-modules.

\subsection{Loewy length and radical, socle series}

Let $V$ be a finite length module over some algebra. 
The Loewy length $ll(V)$ is an invariant of $V$, defined as the length of
the shortest filtration of $V$ with semisimple quotients. Such shortest filtration is not necessarily unique; 
in fact, there are two canonical choices.

The radical of $V$ is defined as the intersection of all maximal proper submodules of $V$. 
The radical filtration, also called the radical series or the upper Loewy series,
\begin{equation}
0 = \rad^l V \subset \rad^{l-1} V \subset \dots \subset\rad^1 V \subset \rad^0 V = V
\end{equation}
is defined inductively by $\rad^{i+1} V = \rad(\rad^i V)$ for $i\ge 0$. 
One can show that the length $l$ of the radical filtration of $V$ is equal to $ll(V)$.

Similarly, the socle of $V$ is defined as the sum of all semisimple submodules of $V$. The socle filtration,
also called the socle series or the lower Loewy series,
\begin{equation}
0 = \soc^0 V \subset \soc^1 V \subset \dots \subset \soc^{l-1} V \subset \soc^l V = V
\end{equation}
is defined inductively by $\soc^{i+1} V/\soc^i V = \soc(V/\soc^i V)$ for $i\ge 0$. Again, the length $l$ of the socle filtration
of $V$ is equal to $ll(V)$.

The semisimple modules $\socl^i V = \soc^i V/\soc^{i-1} V$ and $\radl^i V = \rad^{i-1} V/\rad^i V$
are called the layers of the corresponding series. Moreover, any other filtration 
$0 \subset V_0 \subset \dots\subset V_m = V$ with semisimple quotients satisfies
\begin{equation}
\label{eq:abstract semisimple quotients}
\rad^{m-i} V \subset V_i \subset \soc^i V.
\end{equation}

A module is called rigid, if its socle and radical series coincide; in that case they give the unique shortest
filtration with semisimple quotients, which we simply call the Loewy series.

Examples of rigid modules in the category $\mathcal O$ 
include Verma modules $M_\lambda$ for all $\lambda$, and the big projective modules $P_\lambda$ for anti-dominant $\lambda$.
For more information we refer the reader to \cite{Ir} and references therein.

\section{The regular representation on the big cell}

\subsection{The regular representations}

Let $G$ be the simple complex Lie group, corresponding to $\g$. Denote by $\mathfrak R(G)$ the algebra of regular functions on $G$; it has the structure
of a $G \times G$-module, given by
\begin{equation} \label{eq:group regular actions}
(\varrho_1(g)\psi)(x) = \psi(g^{-1} \, x), \qquad
(\varrho_2(g)\psi)(x) = \psi(x \, g), \qquad
g, x \in G.
\end{equation}

The Lie algebra $\g$ can be defined as the tangent space to the group at the unit element $e \in G$.
To any $\xi \in \g$, we can associate two vector fields on $G$: the left-invariant vectors field 
$\mathcal L_\xi$ such that $\mathcal L_\xi(e) = \xi$, and the right-invariant vector
field $\mathcal R_{-\xi}$ such that $\mathcal R_{-\xi}(e) = -\xi$.
These maps yield two commuting embeddings of $\g$ into the Lie algebra $\operatorname{Vect} (G)$,
and define the two regular $\g$-actions on $\mathfrak R(G)$ by first order differential operators.
Equivalently, they can be defined by the infinitesimal versions 
of \eqref{eq:group regular actions}:
\begin{equation} 
\label{eq:algebra regular actions}
(\varrho_1(\xi) \psi)(x) = \frac d{dt} \psi(e^{-t \, \xi} g ) \biggr|_{t=0}, \qquad
(\varrho_2(\xi) \psi)(x) = \frac d{dt} \psi(g e^{t \, \xi} ) \biggr|_{t=0},
\qquad \xi \in \g, \, x \in G.
\end{equation}

Under both regular actions, the elements of the center $Z(\g)$ act on $\mathfrak R(G)$ by $G\times G$-invariant differential operators; 
for example, the quadratic Casimir operator in $Z(\g)$ corresponds to the Laplace operator on $G$.
The left and right actions of $Z(\g)$ on $\mathfrak R(G)$ differ by an involution $z \mapsto z^*$,
induced by the antipode $x \mapsto x^* = -x$ for $x \in \g$, cf. \eqref{eq:antipode involution}:
\begin{equation}
\varrho_1(z) \varphi = \varrho_2(z^*) \varphi, \qquad z \in Z(\g),\, \varphi \in \mathfrak R(G),
\end{equation}
Equivalently, in terms of the Harish-Chandra isomorphism $Z(\g) \cong S(\h)^W$, it corresponds to the involution 
of $S(\h)^W$ induced by the map $\mathbf h \mapsto -\mathbf h$ for $\mathbf h \in \h$.

Similarly to \eqref{eq:category O blocks}, we define the $\ggbar$-submodules $\mathfrak R(G)_\lambda$ for any $\lambda \in \h^*$ by
\begin{equation}
\label{eq:R(G) block description}
\mathfrak R(G)_\lambda = \left\{\varphi \in \mathfrak R(G) \, \biggr| \, \forall z \in Z(\g) \exists n \in \mathbb N: \ 
(\varrho_1(z^*)-\chi_{\lambda})^n \varphi = (\varrho_2(z)-\chi_{\lambda})^n \varphi = 0 \right\} .
\end{equation}
Since $\mathfrak R(G)$ is locally $G \times G$-finite, it decomposes into the direct sum
\begin{equation}
\label{eq:R(G) block decomposition}
\mathfrak R(G) = \bigoplus_{\lambda \in \mathbf P^+} \mathfrak R(G)_\lambda
\end{equation}
of submodules, corresponding to dominant integral highest weights $\lambda \in \mathbf P^+$. 

Let $T$ and $N_{\pm}$ denote the maximal torus and the unipotent subgroups of $G$, corresponding to $\h$ and $\n_\pm$.
The Gauss decomposition determines the dense open subset $G_0 = N_- \, T \, N_+$ of the group $G$. We refer to $G_0$ as
the big cell of the Gauss decomposition of $G$, and denote by $\mathfrak R(G_0)$ the algebra of regular functions on $G_0$.

The group $G$  does not act on the algebra $\mathfrak R(G_0)$, because $G_0$ is not invariant under left and right shifts by $G$.
Nevertheless, the infinitesimal regular actions \eqref{eq:algebra regular actions} of $\g$ on $\mathfrak R(G_0)$ are still well-defined.
We define the $\ggbar$-submodules $\mathfrak R(G_0)_\lambda$ the same way as in \eqref{eq:R(G) block description},
and obtain the direct sum decomposition analogous to \eqref{eq:R(G) block decomposition}:
\begin{equation}
\label{eq:R(G0) block decomposition}
\mathfrak R(G_0) = \bigoplus_{\lambda \in -\mathbf P^{++}} \mathfrak R(G_0)_\lambda
\end{equation}
The fact that the summation above is over the strictly anti-dominant $\lambda \in -\mathbf P^{++}$
immediately follows from the explicit polynomial realization of $\mathfrak R(G_0)$, which we consider in the next subsection.

\subsection{The polynomial realization of $\mathfrak R(G_0)$}

The algebra $\mathfrak R(G_0)$ can be realized as a certain polynomial algebra, with the
action of $\g$ given by first order differential operators.

It is convenient to order the positive roots $\beta_1,\dots,\beta_m$ according to their
heights,
\begin{equation}
i<j  \ \Rightarrow \ \hgt \beta_i \le \hgt \beta_j, \qquad   i,j = 1,\dots, m = |\Delta_+|.
\end{equation}
and choose the coordinate system $G_0 = N_- \, T \, N_+$ by using the parameterization
\begin{equation}
\label{eq:big cell coordinates}
g(\vec x, \vec y, \vec z) = \exp(x_{\beta_m}\, \mathbf f_{\beta_m}) \dots \exp(x_{\beta_1}\, \mathbf f_{\beta_1}) \,
 z_1^{\mathbf h_1} \dots z_r^{\mathbf h_r} \,
\exp(y_{\beta_1}\, \mathbf e_{\beta_1}) \dots \exp(y_{\beta_m}\, \mathbf e_{\beta_m}),
\end{equation}
where $\vec z = \{z_i\}\in (\C^\times)^r$ are coordinates on $T$, and $\vec x = \{x_\beta\}\in\C^m$ 
(resp. $\vec y = \{y_\beta\}\in\C^m$) are coordinates on $N_-$ (resp. $N_+$).
We also abbreviate $x_i \equiv x_{\alpha_i}$ and $y_i \equiv y_{\alpha_i}$ for $i=1,\dots,r$.
The Gauss decomposition yields the algebra isomorphisms
\begin{equation}
\mathfrak R(G_0) \cong \mathfrak R(N_-) \o  \mathfrak R(T) \o \mathfrak R(N_+) \cong \C[\vec x] \o \C[\mathbf P] \o \C[\vec y],
\end{equation}
where $\mathfrak R(T) = \C[z_1^{\pm1},\dots,z_r^{\pm1}]$ is identified with the group algebra $\C[\mathbf P]$ by
\begin{equation}
z^\lambda \  \equiv \ z_1^{\<\lambda,\mathbf h_1\>} \dots z_r^{\<\lambda,\mathbf h_r\>} \qquad \text{ for any } \lambda \in \mathbf P.
\end{equation}

In general the explicit formulas for the differential
operators describing the action of $\g$ are quite complicated. For our purposes we 
only need a modest amount of information about the structure of this realization.

\begin{prop}
\label{thm:explicit regular actions}
There exist polynomials $p_{i,\beta}, q_{i,\beta}, r_{i,\beta}, s_{i,\beta}$ for $i=1,\dots,r$ and $\beta \in \Delta_+$, 
such that the $\ggbar$-action on $\mathfrak R(G_0)$ is given by
\begin{equation}
\label{eq:explicit left action}
\begin{split}
\varrho_1(\mathbf e_i) &= x_i^2 \, \pd{x_i} - x_i \, z_i \, \pd{z_i} - z_i^{-2} \, \pd{y_i} - \sum_{\beta\in\Delta_+\setminus\{\a_i\}} r_{i,\beta}(\vec x)\, \pd{x_\beta} - z_i^{-2} \sum_{\beta\in\Delta_+\setminus\Pi} s_{i,\beta}(\vec y) \, \pd {y_\beta},\\
\varrho_1(\mathbf h_i) &= - z_i \, \pd{z_i} + 2 x_i \, \pd{x_i} - \sum_{\beta\in \Delta_+\setminus\{\a_i\}} q_{i,\beta}(\vec x)\,\pd{x_\beta},\\
\varrho_1(\mathbf f_i) &= - \pd{x_i} - \sum_{\beta\in \Delta_+\setminus\Pi} p_{i,\beta}(\vec x)\, \pd{y_\beta},\\
\end{split}
\end{equation}

\begin{equation}
\label{eq:explicit right action}
\begin{split}
\varrho_2(\mathbf e_i) &= \pd{y_i} + \sum_{\beta\in \Delta_+\setminus\Pi} p_{i,\beta}(\vec y)\, \pd{y_\beta},\\
\varrho_2(\mathbf h_i) &=  z_i \, \pd{z_i} - 2 y_i \, \pd{y_i} + \sum_{\beta\in \Delta_+\setminus\{\a_i\}} q_{i,\beta}(\vec y)\,\pd{y_\beta},\\
\varrho_2(\mathbf f_i) &= -y_i^2 \, \pd{y_i} + y_i \, z_i \, \pd{z_i} + z_i^{-2} \, \pd{x_i} + \sum_{\beta\in\Delta_+\setminus\{\a_i\}} r_{i,\beta}(\vec y)\, \pd{y_\beta} + 
z_i^{-2} \sum_{\beta\in\Delta_+\setminus\Pi} s_{i,\beta}(\vec x) \, \pd {x_\beta}
\end{split}
\end{equation}

\end{prop}

\begin{proof}
The algorithm for computing these explicit actions is to push the infinitesimal exponentials in 
\eqref{eq:algebra regular actions} inside until they are absorbed by the corresponding factor in \eqref{eq:big cell coordinates}.
To handle the adjustment factors appearing in this commutation process, it suffices to do the computations modulo $t^2$ using the identity 
$$\exp(B) \exp(tA) \equiv  \exp \(t \sum_{j=0}^\infty \frac 1{j!} \underbrace{[B,\dots,[B,[B,A]]\dots]}_{j \text{ commutators}}\) \exp(B) \mod t^2,$$
and similar special cases of the Campbell-Hausdorff formula. Using this approach, one can verify that the right regular action has the form
\eqref{eq:explicit right action} for some $p_{i,\beta}, q_{i,\beta}, r_{i,\beta}, s_{i,\beta}$; we skip the technical details.

It is easy to see that the transposition anti-automorphism $g \mapsto g^{\top}$ interchanges the $\vec x$ and $\vec y$ variables, 
i.e. $g(\vec x,\vec y, \vec z)^\top = g(\vec y, \vec x, \vec z)$, and transforms the right regular action $\varrho_2$ into
the negative left regular action $-\varrho_1$.
This implies that the left regular action \eqref{eq:explicit left action}same polynomials is described in terms of the same polynomials 
$p_{i,\beta}, q_{i,\beta}, r_{i,\beta}, s_{i,\beta}$.
\end{proof}

In other words, Proposition \ref{thm:explicit regular actions} states that for each simple root $\a_i$ the action on $\mathfrak R(G_0)$ of the subalgebra 
$\mathfrak{sl}_{\a_i}(2,\C) = \C\mathbf e_i \oplus \C\mathbf h_i \oplus\C\mathbf f_i$ is given by the standard $\mathfrak{sl}(2,\C)$ action
(see e.g. \cite{FS}), plus additional terms determined by $p_{i,\beta}, q_{i,\beta}, r_{i,\beta}, s_{i,\beta}$. Other properties of these polynomials 
- for example, they are homogeneous in the appropriate sense - are not crucial for our further analysis.

\begin{prop}
For any $\lambda \in - \mathbf P^{++}$ we have
\begin{equation}
\label{eq:R(G_0) character}
\ch \mathfrak R(G_0)_\lambda = \sum_{w\in W^\lambda} \ch \( M_{w \cdot \lambda}^* \o M_{w \cdot \lambda} \).
\end{equation}
Equivalently, for any $x,y \in W^\lambda$ we have
\begin{equation}
\label{eq:R(G_0) JH multiplicities}
[\mathfrak R(G_0)_\lambda: L_{x \cdot\lambda} \o L_{y \cdot\lambda}^*] = [P_{y\cdot\lambda} : L_{x \cdot\lambda}].
\end{equation}
\end{prop}
\begin{proof}
Consider the collection  $\{\mathfrak R^{\le \lambda}\}_{\lambda\in\mathbf P}$ of subspaces of $\mathfrak R(G_0)$, defined by
\begin{equation}
\mathfrak R^{\le \lambda} = \bigoplus_{\mu\le\lambda} \C[\vec x] \o \C z^\mu \o \C[\vec y]
\end{equation}
where $\mu \le \lambda$ is equivalent to $\lambda-\mu \in \mathbf P^+$. It is clear from \eqref{eq:explicit left action}, \eqref{eq:explicit right action} that
$\mathfrak R^{\le \lambda}$ is a $\ggbar$-submodule of $\mathfrak R(G_0)$ for any $\lambda \in \mathbf P$. 
Similarly, one defines the $\ggbar$-submodules $\{\mathfrak R^{< \lambda}\}_{\lambda\in\mathbf P}$ of $\mathfrak R(G_0)$.

The induced $\ggbar$-actions in the quotients $\mathfrak R^{\le\lambda}/\mathfrak R^{<\lambda}$ are given by
\begin{align}
\label{eq:flag variety action 1}
\begin{split}
\varrho_1^{(\lambda)}(\mathbf e_i) &= x_i^2 \, \pd{x_i} - x_i \, z_i \, \pd{z_i} - 
\sum_{\beta\in\Delta_+\setminus\{\a_i\}} r_{i,\beta}(\vec x)\, \pd{x_\beta} ,\\
\varrho_1^{(\lambda)}(\mathbf h_i) &= - z_i \, \pd{z_i} + 2 x_i \, \pd{x_i} - \sum_{\beta\in \Delta_+\setminus\{\a_i\}} q_{i,\beta}(\vec x)\,\pd{x_\beta},\\
\varrho_1^{(\lambda)}(\mathbf f_i) &= - \pd{x_i} - \sum_{\beta\in \Delta_+\setminus\Pi} p_{i,\beta}(\vec x)\, \pd{y_\beta},\\
\end{split}
\end{align}
\begin{align}
\label{eq:flag variety action 2}
\begin{split}
\varrho_2^{(\lambda)}(\mathbf e_i) &= \pd{y_i} + \sum_{\beta\in \Delta_+\setminus\Pi} p_{i,\beta}(\vec y)\, \pd{y_\beta},\\
\varrho_2^{(\lambda)}(\mathbf h_i) &=  z_i \, \pd{z_i} - 2 y_i \, \pd{y_i} + \sum_{\beta\in \Delta_+\setminus\{\a_i\}} q_{i,\beta}(\vec y)\,\pd{y_\beta},\\
\varrho_2^{(\lambda)}(\mathbf f_i) &= -y_i^2 \, \pd{y_i} + y_i \, z_i \, \pd{z_i} + \sum_{\beta\in\Delta_+\setminus\{\a_i\}} r_{i,\beta}(\vec y)\, \pd{y_\beta} 
\end{split}
\end{align}
The $\g$-action \eqref{eq:flag variety action 1} is homogeneous in $\vec z$ and independent of $\vec x$, and therefore defines an action of $\g$ on the space $\C z^\lambda  \otimes \C[\vec x]$, which is identified with the functions (or sections of line bundles) on the big cell of the flag variety. The corresponding $\g$-module can be shown to be is isomorphic to $M_\lambda^*$. Similarly, the action \eqref{eq:flag variety action 2} in $\C z^\lambda \o \C[\vec y]$
yields the contragredient Verma module $M_\lambda^c$.
Therefore, we have a $\ggbar$-module isomorphism
\begin{equation}
\mathfrak R^{\le\lambda}/\mathfrak R^{<\lambda} = \C[\vec y] \o \C z^\lambda \o \C[\vec x] \cong M_\lambda^* \o M_\lambda^c.
\end{equation}
Collecting the constituents associated with the central character $\chi_\lambda$, i.e. arising from the quotients $\mathfrak R^{\le\mu}/\mathfrak R^{<\mu}$ with $\mu \in W^\lambda \cdot \lambda$, we obtain
\begin{equation}
\ch \mathfrak R(G_0)_\lambda = \sum_{w\in W^\lambda} \ch \( M_{w \cdot \lambda}^* \o M_{w \cdot \lambda}^c \) =
\sum_{w\in W^\lambda} \ch \( M_{w \cdot \lambda}^* \o M_{w \cdot \lambda} \),
\end{equation}
establishing \eqref{eq:R(G_0) character}. Similarly, we have
\begin{multline*}
[\mathfrak R(G_0)_\lambda: L_{x \cdot\lambda} \o L_{y \cdot\lambda}^*] = 
\sum_{w\in W^\lambda} [M_{w \cdot \lambda}^* \o M_{w \cdot \lambda}^c : L_{x\cdot \lambda}^* \o L_{y\cdot \lambda}] = \\
= \sum_{w\in W^\lambda} [M_{w \cdot \lambda}^*: L_{x\cdot \lambda}^*] [M_{w \cdot \lambda}^c : L_{y\cdot \lambda}] = 
\sum_{w\in W^\lambda} [M_{w \cdot \lambda}  : L_{x\cdot \lambda}] [P_{y \cdot \lambda} : M_{w\cdot \lambda}] = [P_{y \cdot \lambda} : L_{x\cdot \lambda}].
\end{multline*}
\end{proof}

\section{Whittaker vectors in $\mathfrak R(G_0)$ and the big projective modules}

\subsection{Whittaker vectors and Whittaker functor}

Let $\eta = (\eta_1,\dots, \eta_r) \in \C^r$. There exists a unique character
$\boldsymbol\eta^+: \mathcal U(\n_+) \to \C$, such that $\boldsymbol\eta^+(\mathbf e_i) = \eta_i$ for $i=1,\dots,r$.
In order for $\boldsymbol \eta^+$ to be an algebra homomorphism we must put $\boldsymbol\eta^+(\mathbf e_\beta) = 0$ for $\beta \notin \Pi$.

The Whittaker functor $\wh_\eta^+$ associates to any $\g$-module $V$ the subspace
\begin{equation}
\wh_\eta^+ (V) = \left\{v \in V \, \bigr| \, \ker \boldsymbol\eta^+ \cdot v = 0 \right\}.
\end{equation}
Elements $v \in \wh_\eta^+ (V)$ are called the Whittaker vectors of $V$ with respect to $\eta$.

One can show that two Whittaker functors $\wh_\eta^+$ and $\wh_{\eta'}^+$ are isomorphic
if and only if the corresponding characters vanish on the same Chevalley generators, 
i.e. $\eta_i = 0 \Leftrightarrow \eta'_i = 0$.
Thus, there are $2^r$ nonequivalent Whittaker functors. 
The two extreme cases correspond to the trivial character with to all $\eta_i=0$, 
and to the nonsingular characters with all $\eta_i\ne0$.

In the case of trivial $\eta=0$, the functor $\wh_0^+$ gives the subspace of singular vectors in $V$:
\begin{equation}
\wh_0^+ (V) = \operatorname{Sing}^+ (V) = \left\{v \in V \, \bigr| \, \n_+ \cdot v = 0 \right\},
\end{equation}
which are important in the study of modules in category $\mathcal O$.
However, for nontrivial $\eta$ we have $\wh_\eta^+ V = 0$ for all $V \in \mathcal O$,
which is easily proved using the weight decomposition of $V$.

Therefore, we should consider a modification $\Wh_\eta^+$ of the Whittaker functor $\wh_\eta^+$.
For any $V = \bigoplus_\mu V[\mu]$ introduce its completion $\bar V = \prod_\mu V[\mu]$.
We define
\begin{equation}
\Wh_\eta^+(V) = \wh_\eta^+( \bar V ) = \left\{v \in \bar V \, \bigr| \, \ker \boldsymbol\eta^+ \cdot v = 0 \right\}.
\end{equation}

The center $Z(\g)$ of the universal enveloping algebra $\mathcal U(\g)$ acts on the space of Whittaker vectors in $V$.
Therefore, for each $\eta$ we get a functor
\begin{equation}
\Wh_\eta^+: \g\text{-mod} \to Z(\g)\text{-mod}, \qquad V \mapsto \Wh_\eta^+ (V).
\end{equation}
We will be particularly interested in the restriction of the functor $\Wh_\eta^+$ to the category $\mathcal O$.

Similarly, one defines versions $\wh_\eta^-,\Wh_\eta^-: \g\text{-mod} \to Z(\g)\text{-mod}$ of the Whittaker functors, 
corresponding to the character 
$\boldsymbol\eta^-: \mathcal U(\n_-) \to \C$ defined by $\boldsymbol\eta^- (\mathbf f_i) = \eta_i$. 
The functor $\Wh_\eta^-$ plays the same role for the category $\mathcal O^*$ as $\Wh_\eta^+$ plays for $\mathcal O$.

If $V$ is equipped with two commuting $\g$-actions, then the space of Whittaker vectors associated with 
one copy of $\g$ retains the $\g$-module structure with respect to the other copy of $\g$. In the next section
we study the functors
\begin{alignat*}{5}
\Wh_\eta^- \o 1: \ & \ggbar\text{-mod} &\ \to & \ Z(\g) \oplus \g\text{-mod},\\
1 \o \Wh_\eta^+: \ & \ggbar\text{-mod} &\ \to & \ \g \oplus Z(\g)\text{-mod},\\
\Wh_\eta^- \o \Wh_{\eta'}^+: \ & \ggbar\text{-mod} &\ \to & \ Z(\g) \oplus Z(\g)\text{-mod},
\end{alignat*}
applied to the regular representation on the big cell.

\subsection {Polynomial realization of $(\Wh_\eta^- \o 1)(\mathfrak R(G_0))$}

We show that the space of Whittaker vectors with respect to the left regular action admits a polynomial realization,
where the right regular action of $\g$ is very similar to the $\g$-actions on functions on the flag variety.

\begin{prop}
The $\g$-module $(\Wh_\eta^- \o 1)(\mathfrak R(G_0))$ can be realized in the polynomial algebra $\C[\mathbf P] \o \C[\vec y]$, with
the $\g$-action given by
\begin{equation}
\label{eq:explicit projective action}
\begin{split}
\varrho(\mathbf e_i) &= \pd{y_i} + \sum_{\beta\in \Delta_+\setminus\Pi} p_{i,\beta}(\vec y)\, \pd{y_\beta},\\
\varrho(\mathbf h_i) &=  z_i \, \pd{z_i} - 2 y_i \pd{y_i} + \sum_{\beta\in \Delta_+\setminus\{\a_i\}} q_{i,\beta}(\vec y)\,\pd{y_\beta},\\
\varrho(\mathbf f_i) &= -y_i^2 \pd{y_i} + y_i \, z_i \pd{z_i}  + 
\sum_{\beta\in\Delta_+\setminus\{\a_i\}} r_{i,\beta}(\vec y)\, \pd{y_\beta} + \eta_i \, z_i^{-2}
\end{split},
\end{equation}
where $p_{i,\beta}, q_{i,\beta}, r_{i,\beta}$ are the same as in Proposition \ref{thm:explicit regular actions}.
\end{prop}

\begin{proof}
The desired Whittaker vectors in $\mathfrak R(G_0)$ are represented by functions satisfying
\begin{equation}
\label{eq:Whittaker DE}
\(\pd{x_i} + \sum_{\beta \in \Delta^+} p_{i,\beta}(\vec x) \pd{x_\beta} + \eta_i \) \psi(\vec x, \vec y, \vec z) = 0, \qquad i=1,\dots,r
\end{equation}
These differential equations are independent of $\vec y, \vec z$, and thus the space of solutions is
spanned by functions of the form $\psi(\vec x, \vec y, \vec z) = \phi(\vec y, \vec z) \ \tau(\vec x)$
where $\tau(\vec x)$ satisfies the equations \eqref{eq:Whittaker DE}, and represents a Whittaker vector in $\C[\vec x] \cong M_0^*$. It is known that 
for any $\mu \in \h^*$ and $\eta\in \C^r$ we have $\dim \Wh_\eta^-(M_\mu^*) = 1$, see e.g. \cite{Ba}.
Therefore, there exists a unique up to proportionality solution $\tau(\vec x) \in \C[[x]]$ of \eqref{eq:Whittaker DE}, and a
direct verification shows that $\tau(\vec x) = \exp \( \sum_{i=1}^r \eta_i \, x_{\a_i} \)$.

We conclude that the Whittaker vectors in $\mathfrak R(G_0)$ are identified with functions
\begin{equation}
\label{eq:left Whittaker functions}
\psi(\vec x, \vec y, \vec z) = \phi(\vec y, \vec z) \ \exp \( \sum_{i=1}^r \eta_i \, x_{\a_i} \),
\end{equation}
and it is straightforward to check that the specialization of the formulas \eqref{eq:explicit right action}
to the subspace of functions \eqref{eq:left Whittaker functions} results in the action
\eqref{eq:explicit projective action} on the space $\C[\mathbf P] \o \C[\vec y]$.
\end{proof}

The polynomial realization gives us the following important information on the size of the Whittaker vector spaces.

\begin{prop}
For any $\eta \in \C^r$ we have
\begin{equation}
\label{eq:projective block character}
\ch (\Wh_{\eta}^- \o 1) (\mathfrak R(G_0)_\lambda)  = \sum_{w \in W^\lambda} \ch M_{w\cdot\lambda}.
\end{equation}
Equivalently, for any $w \in W^\lambda$ we have
\begin{equation}
\label{eq:projective block composition}
[(\Wh_{\eta}^- \o 1) (\mathfrak R(G_0)_\lambda) : L_{x\cdot\lambda}] = [P_\lambda : L_{x\cdot\lambda}].
\end{equation}
\end{prop}

\begin{proof}
Consider the collection of submodules $\mathfrak P^{\le\lambda}$, defined by
\begin{equation}
\mathfrak P^{\le \lambda} = \bigoplus_{\mu\le\lambda} \C z^\mu  \o \C[\vec y].
\end{equation}
The corresponding quotients $\mathfrak P^{\le\lambda}/\mathfrak P^{<\lambda}$ are given by
\eqref{eq:flag variety action 1}, and therefore
\begin{equation}
\label{eq:projective realization quotient} 
\mathfrak P^{\le\lambda}/\mathfrak P^{<\lambda} \cong M_\lambda^c.
\end{equation}

Collecting the constituents from $\mathcal O_\lambda$, we conclude that
\begin{equation}
\ch (\Wh_{\eta}^- \o 1) (\mathfrak R(G_0)_\lambda)  = \sum_{w \in W^\lambda} \ch M_{w\cdot\lambda}^c = 
\sum_{w \in W^\lambda} \ch M_{w\cdot\lambda}.
\end{equation}
Similarly, using the fact that for every $w \in W^\lambda$ we have $[P_\lambda:M_{w\cdot\lambda}] = 1$, we obtain
\begin{equation}
[(\Wh_{\eta}^- \o 1) (\mathfrak R(G_0)_\lambda) : L_{x\cdot\lambda}] = 
\sum_{w\in W^\lambda} [M_{w\cdot\lambda}^c: L_{x\cdot\lambda}] = [P_\lambda : L_{x\cdot\lambda}] .
\end{equation}
\end{proof}

\subsection{Big projective modules and Borel-Weil realization}

The Borel-Weil realization of finite-dimensional simple $\g$-modules $L_\lambda$ with dominant integral highest weights $\lambda\in\mathbf P^+$,
given by \eqref{eq:Borel-Weil realization}, can be reformulated in the following form.

\begin{thm}
For any $\lambda \in \mathbf P^+$, there is an isomorphism
\begin{equation}
L_\lambda \cong (\operatorname{Sing}^- \o 1) (\mathfrak R(G)_\lambda) = \left\{\varphi \in \mathfrak R(G)_\lambda \, \biggr| \, \varrho_1(\n_-) \varphi = 0 \right\}.
\end{equation}
\end{thm}

Replacing the functor $\operatorname{Sing}^-$ by its generalization $\Wh_\eta^-$ and applying it to the regular representation on the big cell,
we obtain a similar realization of the big projective modules in the category $\mathcal O$.

\begin{thm}
\label{thm:block Borel-Weil}
Let $\eta \in \C^r$ be nonsingular. Then for any $\lambda \in -\mathbf P^{++}$ there is an isomorphism
\begin{equation}
\label{eq:block Borel-Weil}
(\Wh_{\eta}^- \o 1)(\mathfrak R(G_0)_\lambda) \cong P_{\lambda}.
\end{equation}
\end{thm}

\begin{proof}
Set for brevity $V = (\Wh_{\eta}^- \o 1)(\mathfrak R(G_0)_\lambda)$. The proof is based on two properties of $V$. 

First, the module $V$ has a quotient, isomorphic to the contragredient Verma module $M_{w_\circ\cdot\lambda}^c$,
which immediately follows from \eqref{eq:projective realization quotient}.

Second, $V$ has a submodule, isomorphic to the Verma module $M_{w_\circ\cdot\lambda}$.
Indeed, consider the vector $v \in V$ represented by the polynomial $z^{w_\circ\cdot\lambda}$; it is annihilated by $\n_+$ and has
dominant highest weight $w_\circ\cdot\lambda$. The submodule generated by $v$ is {\it a priori} isomorphic to
a quotient of $M_{w_\circ\cdot\lambda}$; we show that in fact it coincides with $M_{w_\circ\cdot\lambda}$.

The socle of the Verma module $M_{w_\circ\cdot\lambda}$ is generated by a singular vector $v_{sing} = q(\lambda) v_0$, where $v_0$ is the highest weight
vector generating $M_{w_\circ\cdot\lambda}$, and $q(\lambda)$ is some element in $\mathcal U(\n_-)$ of weight $w_\circ\cdot\lambda -\lambda$.
Let $\{n_i\}_{i=1}^r$ be the nonnegative integers such that $\sum_{i=1}^r n_i \a_i = w_\circ\cdot\lambda - \lambda$. Consider the PBW basis
of $\mathcal U(\n_-)$, associated with the ordering $\{\mathbf f_{\beta_1}, \dots, \mathbf f_{\beta_m} \}$. Then with respect to this basis
we have an expansion
$$q(\lambda) = c \ \mathbf f_1^{n_1} \dots \mathbf f_r^{n_r} + \text{ terms with nonsimple root vectors},$$
with some nonzero coefficient $c$, see e.g. \cite{Ba}. It is clear that in the polynomial realization of $V$ the action of $\mathbf f_\beta$ contains the 
multiplication by $z^{-\beta}$ if and only if $\beta$ is a simple root. Therefore, the coefficient before $z^\lambda$
in the expansion of $q(\lambda) v$  is equal to $c \,  \eta_1^{m_1} \dots \eta_r^{m_r} \ne 0$, and thus 
$q(\lambda) v \ne 0$.
Hence $v$ generates a submodule of $V$ isomorphic to $M_{w_\circ\cdot\lambda}$.

We now return to the main proof. It follows from \eqref{eq:projective block composition} that the composition series for
$V$ contains exactly one simple module $L_{w_\circ\cdot\lambda}$ with the dominant highest weight $w_\circ\cdot\lambda$. 
This unique constituent appears in the socle of $M_\lambda^c$
and in the top layer of $M_\lambda$. It is known that both $M_\lambda$ and $M_\lambda^c$ are rigid and have Loewy length $l_\lambda+1$,
and therefore $V$ has Loewy length at least $2l_\lambda+1$.

On the other hand, modules in the category $\mathcal O_\lambda$ have Loewy Length at most $2l_\lambda+1$, so we have $ll(V) = 2l_\lambda+1$.
Moreover, the only indecomposable module of Loewy length $2l_\lambda+1$ is the big projective module $P_\lambda$. Therefore, $V$ must contain $P_\lambda$ as a direct summand, and comparing the characters \eqref{eq:big projective character} and \eqref{eq:projective block character} 
we see that in fact $V \cong P_\lambda$.
\end{proof}

Similar arguments show that for any $\lambda \in  - \mathbf P^{++}$ there is an isomorphism
\begin{equation}
(1 \o \Wh_{\eta}^+)(\mathfrak R(G_0)_\lambda) \cong P_\lambda^*.
\end{equation}

\subsection{Whittaker functions on $G_0$}

Let $\eta,\eta' \in \C^r$. The associated Whittaker functions on the group $G$ are defined as 
elements of $(\Wh_\eta^- \o \Wh_{\eta'}^+) (\mathfrak R(G_0))$. When $G = SL(2,\C)$,
they are directly related to the Whittaker functions $W_{k,m}(z)$, satisfying the Whittaker
differential equation \cite{WW}. For Lie groups of other type the Whittaker functions
were used by Kazhdan and Kostant to prove the integrability of the quantum Toda system:
the restriction of the Laplacian to the subspace of Whittaker functions coincides
with the Toda Hamiltonian, and higher Casimir operators yield the quantum integrals of motion.

The decomposition \eqref{eq:R(G0) block decomposition} induces the $Z(\g) \o Z(\g)$-module decomposition
\begin{equation}
(\Wh_\eta^- \o \Wh_{\eta'}^+) (\mathfrak R(G_0)) = \bigoplus_{\lambda \in - \mathbf P^{++}}(\Wh_\eta^- \o \Wh_{\eta'}^+) (\mathfrak R(G_0)_\lambda ),
\end{equation}
The subspaces $(\Wh_\eta^- \o \Wh_{\eta'}^+) (\mathfrak R(G_0)_\lambda )$ are spanned by the generalized eigenfunctions
for the quantum Toda system, corresponding to the central character $\chi_\lambda$.

For any $\lambda \in - \mathbf P^{++}$ denote $\mathcal C_\lambda = \End_{\mathcal O} (P_\lambda)$.
The algebra $\mathcal C_\lambda$ was studied in \cite{So,Be}, and can be interpreted as the cohomology
ring of a suitable complex flag variety. In particular, $\mathcal C_\lambda$ is commutative.

\begin{prop}
\label{thm:double Whittaker}
Let $\eta,\eta' \in \C^r$ be nonsingular. Then
\begin{equation}
(\Wh_\eta^- \o \Wh_{\eta'}^+) (\mathfrak R(G_0)_\lambda ) \cong \mathcal C_\lambda.
\end{equation}
\end{prop}

\begin{proof}
Using Theorem \ref{thm:block Borel-Weil}, we get
\begin{equation}
(\Wh_\eta^- \o \Wh_{\eta'}^+) (\mathfrak R(G_0)_\lambda ) = (1 \o \Wh_{\eta'}^+) (\Wh_\eta^- \o 1) (\mathfrak R(G_0)_\lambda ) \cong
\Wh_{\eta'}^+ (P_\lambda).
\end{equation}
It was shown in \cite{Ba} that for nonsingular $\eta$ the restriction of the Whittaker functor $\Wh_\eta^+$ is isomorphic to Soergel's functor
\begin{equation}
\mathbb V: \mathcal O_\lambda \to \mathcal C_\lambda\text{-mod}, \qquad
\mathbb V(M) = \Hom_{\mathcal O}(P_\lambda,M),
\end{equation}
when we regard $\mathbb V$ as a functor from $\mathcal O_\lambda$ to $Z(\g)$-modules via the surjection $Z(\g) \to \mathcal C_\lambda$, arising from the
action of $\mathcal U(\g)$ on $P_\lambda$ (see \cite{So}).
Therefore, we have 
\begin{equation}
\Wh_{\eta'}^+ (P_\lambda) \cong \Hom_{\mathcal O}(P_\lambda,P_\lambda) = \mathcal C_\lambda.
\end{equation}
\end{proof}

Using the equivalence of Whittaker and Soergel's functors, we reformulate Theorem \ref{thm:block Borel-Weil}. 
\begin{cor}
\label{thm:algebraic Whittaker}
For any $\lambda \in \mathbf P^{++}$, there are $\g$-module isomorphisms
\begin{equation}
\label{eq:algebraic Whittaker}
\Hom_{\varrho_1}(P_\lambda^*,\mathfrak R(G_0)) \cong P_\lambda, \qquad
\Hom_{\varrho_2}(P_\lambda,\mathfrak R(G_0)) \cong P_\lambda^*,
\end{equation}
where we regard $\mathfrak R(G_0)$ as a $\g$-module
with respect to the indicated regular action $\varrho_1$ or $\varrho_2$, and $\g$ acts on the space of morphisms via the
other regular action.
\end{cor}

Another consequence of Proposition \ref{thm:double Whittaker} is that $\dim (\Wh_\eta^- \o \Wh_{\eta'}^+) (\mathfrak R(G_0)_\lambda ) = |W^\lambda|$.
Thus, the space of generalized eigenfunctions of the quantum Toda system has dimension $|W^\lambda|$, and
has a natural grading, inherited from the corresponding cohomology ring.
It would be interesting to interpret this grading in the context of quantum completely integrable systems.

\section{Matrix elements of modules in the category $\mathcal O$.}

\subsection{The algebra $\mathcal U(\g)^*$ and subspaces of matrix elements}

Let $\mathcal U(\g)^*$ denote the dual space of the enveloping algebra $\mathcal U(\g)$;
it has a natural $\ggbar$-module structure, defined by
\begin{equation}
(\varrho_1(\xi)\psi)(x) = - \psi(\xi \, x), \qquad 
(\varrho_2(\xi)\psi)(x) = \psi(x \, \xi), \qquad \xi \in \g, \ x \in \mathcal U(\g).
\end{equation}

Let $V$ be a $\g$-module. For any $ v \in V, \, v^* \in V^*$ we define the ``matrix element'' functional
\begin{equation}
\label{eq:matrix elements map}
\Phi_{v^* \o v} (x) = \<v^*, x\, v\>, \qquad x \in \mathcal U(\g),
\end{equation}
and extend it by linearity to the map
\begin{equation}
\Phi: V^* \o V \to \mathcal U(\g)^*.
\end{equation}
We denote by $\mathbb M(V)$ the subspace $\Phi(V^* \o V) \subset \mathcal U(\g)^*$, spanned by the matrix elements
of $V$. It is easy to check that the map $\Phi$ is a $\ggbar$-homomorphism, and thus $\mathbb M(V)$ is a $\ggbar$-submodule of $\mathcal U(\g)^*$.
We  use the same notation $\Phi$ for all $\g$-modules $V$, i.e. regard $\Phi$ as the ``universal'' matrix elements map.

Dualizing the comultiplication in the Hopf algebra $\mathcal U(\g)$, we obtain a commutative associative product in
$\mathcal U(\g)^*$. The commutative algebra $\mathcal U(\g)^*$ is isomorphic to the algebra of the
formal power series in $n=\dim \g$ variables, see e.g. \cite{Di}. 

We consider two smaller subalgebras of $\mathcal U(\g)^*$.
The Hopf dual $\mathcal U(\g)_{Hopf}^*$ is defined by
\begin{equation*}
\mathcal U(\g)_{Hopf}^* = \left\{ \varphi \in \mathcal U(\g)^* \, \biggr| \,
\ker \varphi \text{ contains a two-sided ideal } J\subset\mathcal U(\g) \text { of finite codimension} \right\}.
\end{equation*}
One  immediately checks that $\mathcal U(\g)_{Hopf}^*$ is a subalgebra and a $\ggbar$-submodule of $\mathcal U(\g)^*$.
Moreover, the $\ggbar$ action on $\mathcal U(\g)_{Hopf}^*$ is locally finite, since for any $\varphi \in \mathcal U(\g)_{Hopf}^*$
all elements in $(\mathcal U(\g) \o \mathcal U(\g))(\varphi)$ also annihilate the corresponding ideal $J$, and $\dim \Ann(J) < \infty$.
This is the exact analogue of the locally finite $G \times G$ action on $\mathfrak R(G)$; in fact one has the algebra and $\ggbar$-module
isomorphism $\mathcal U(\g)_{Hopf}^* \cong \mathfrak R(G).$ In other words, $\mathcal U(\g)_{Hopf}^*$ gives the Lie algebraic model
of the regular representation $\mathfrak R(G)$, and in particular is spanned by the matrix elements of simple finite-dimensional
$\g$-modules $L_\lambda$ for $\lambda \in \mathbf P^+$.

The matrix elements of modules in the category $\mathcal O$ span a larger subspace of $\mathcal U(\g)^*$. 
Consider the triangular decomposition $\mathcal U(\g) = \mathcal U(\n_-) \o \mathcal U(\h) \o \mathcal U(\n_+)$, 
and let $\mathcal U(\n_\pm)^*$ be the restricted duals of the corresponding enveloping algebras with respect to the principal gradation.
We also consider the subalgebra of $\mathcal U(\h)^*$, spanned by functionals $e^\mu$, defined by
$$\<e^\mu, \mathbf h_1^{k_1}\dots\mathbf h_r^{k_r}\> = \<\mu,\mathbf h_1\>^{k_1} \dots  \<\mu,\mathbf h_r\>^{k_r},
\qquad \mu \in \h^*, \, k_1,\dots,k_r \in \Z_{\ge0}.$$
We note that under the identification of $\mathcal U(\h)^*$ with the algebra of formal power series, the functionals
$e^\mu$ correspond to the exponential functions, and hence naturally correspond to the elements of the group algebra $\C[\h^*]$.
We define
\begin{equation}
\mathcal U(\g)_{\mathcal O}^* = \mathcal U(\n_-)^* \o \C[\h^*] \o \mathcal U(\n_+)^*.
\end{equation}
The subspace $\mathcal U(\g)_{\mathcal O}^*$ is also a subalgebra and a $\ggbar$-submodule of $\mathcal U(\g)^*$.
It follows from the definition that the $\ggbar$-action on it is locally $\n_- \oplus \n_+$-finite, and that
$\mathcal U(\g)_{\mathcal O}^*$ is $\h$-diagonalizable with respect to both left and right actions.

\begin{prop}
\label{thm:matrix elements in O}
The subalgebra $\mathcal U(\g)^*_{\mathcal O}$ coincides with the subspace $\mathbb M_{\mathcal O}$, defined by
\begin{equation}
\mathbb M_{\mathcal O} = \sum_{V \in \mathcal O} \mathbb M(V) \subset \mathcal U(\g)^*.
\end{equation}
Equivalently, $\mathcal U(\g)^*_{\mathcal O}$ is spanned by the matrix elements of modules in category $\mathcal O$.
\end{prop}

\begin{proof}
Suppose $V \in \mathcal O$, and let $v \in V, \, v^* \in V^*$ be $\h$-homogeneous.
For any $\h$-homogeneous $x_\pm \in \mathcal U(\n_\pm)$ and any $x_0 = \mathbf h_1^{k_1}\dots\mathbf h_r^{k_r} \in \mathcal U(\h)$ we have
\begin{equation}
\Phi_{v^* \o v}(x_- x_0 x_+) = \<v^*, x_- x_0 x_+ v\>  = \<\mu,\mathbf h_1\>^{k_1} \dots  \<\mu,\mathbf h_r\>^{k_r} \ \<x_-^* v^*, x_+v\>,
\end{equation}
where $\mu\in\h^*$ is the weight of the homogeneous vector $x_+ v$. This implies that $\Phi_{v^* \o v} \in \mathcal U(\g)^*_{\mathcal O}$,
and therefore $\mathbb M(V) \subset \mathcal U(\g)^*_{\mathcal O}$.

We now prove the opposite inclusion. Let $\varphi \in \mathcal U(\g)^*_{\mathcal O}$, and consider the subspace 
\begin{equation}
\label{eq:subspace generated by phi}
V = (1 \o \mathcal U(\g)) \, \varphi \subset \mathcal U(\g)^*
\end{equation}
With respect to the right regular action $V$ is an $\h$-diagonalizable $\g$-module, generated by a single $\n_+$-nilpotent element $\varphi$, 
and hence $V \in \mathcal O$.
Let $\varphi^*$ denote the restriction of $1 \in \mathcal U(\g)$ to the subspace $V$, i.e. $\<\varphi^*,\psi\> = \psi(1)$; then $\varphi^* \in V^*$ and
\begin{equation}
\Phi_{\varphi^* \o \varphi}(x) = \<\varphi^*, x \, \varphi\> = (x \varphi)(1) = \varphi(x),
\end{equation}
which means that $\varphi \in \mathbb M(V)$, and completes the proof.
\end{proof}

The subalgebra $\mathcal U(\g)^*_{\mathcal O}$ decomposes into a direct sum of $\ggbar$-submodules, corresponding to the
central characters of $\g$. Let $\mathbb M_\lambda$ be subspace of matrix elements of modules in $\mathcal O_\lambda$,
\begin{equation}
\mathbb M_\lambda = \sum_{V \in \mathcal O_\lambda} \mathbb M(V), \qquad \lambda \in \h^*.
\end{equation}

Another important subalgebra of $\mathcal U(\g)^*$ is the subspace
\begin{equation}
\mathcal U(\g)^*_{\mathcal O_{int}} = \mathcal U(\n_-)^* \o \C[\mathbf P] \o \mathcal U(\n_+)^* \cong
 \bigoplus_{\lambda \in -\mathbf P^{++}} \mathbb M_\lambda, 
\end{equation}
corresponding to the matrix elements of module $V \in O$ with integral weights. This subspace gives the Lie algebraic
model of \eqref{eq:R(G0) block decomposition}, and in fact is isomorphic to $\mathfrak R(G_0)$ as an algebra and as a $\ggbar$-module, see 
Theorem \ref{thm:projective Peter-Weyl}.

\subsection{Rigidity and Loewy series of $\mathbb M_\lambda$}

For any $\lambda \in - \mathbf P^{++}$ the $\ggbar$-module $\mathbb M_\lambda$ admits an increasing filtration 
$0 = \mathbb M_\lambda^{(0)} \subset \mathbb M_\lambda^{(1)} \subset \dots\subset 
\mathbb M_\lambda^{(2l_\lambda)} \subset \mathbb M_\lambda^{(2l_\lambda+1)} = \mathbb M_\lambda$, defined by
\begin{equation}
\label{eq:matrix elements filtration}
\mathbb M_\lambda^{(k)} = \sum_{\substack{V \in \mathcal O_\lambda,\\ ll(V)\le k}}  \mathbb M(V), \qquad k = 1,\dots,2l(w_\lambda)+1.
\end{equation}
Since the Loewy length of any $V \in \mathcal O_\lambda$ is at most $2l_\lambda+1$, we have $\mathbb M_\lambda^{(2l_\lambda+1)} = \mathbb M_\lambda$.

\begin{prop}
\label{thm:M(lambda) socle series}
For any $\lambda \in -\mathbf P^{++}$ the socle filtration of $\mathbb M_\lambda$ coincides with \eqref{eq:matrix elements filtration}.
\end{prop}
\begin{proof}
Consider the layers $\overline{\mathbb M}_\lambda^{(k)} = \mathbb M_\lambda^{(k)}/ \mathbb M_\lambda^{(k-1)}$ of this filtration. 
Let $\varphi \in \mathbb M_\lambda^{(k)}$, and let $V \in \mathcal O_\lambda$ be such that $ll(V)\le k$ and $\varphi \in \mathbb M(V)$.
Since $ll(\rad V) = ll(V)-1 \le k-1$, we have 
\begin{equation*}
\Phi(V^* \o \rad V) = \Phi( (\rad V)^* \o \rad V) \subset \mathbb M^{(k-1)},
\end{equation*}
and similarly $\Phi(\rad V^* \o V) \subset \mathbb M^{(k-1)}$. Therefore, we can complete the commutative diagram
\begin{equation}
\label{eq:diagram factors through}
\begin{gathered}
\xymatrix{
V^* \o V \ar[rrr]^\Phi \ar[d]_{\pi\o\pi} &&& \mathbb M_\lambda^{(k)} \ar[d]^\pi \\
(V^*/\rad V^*) \o (V /\rad V) \ar@{.>}[rrr] &&&   \mathbb M_\lambda^{(k)}/ \mathbb M_\lambda^{(k-1)}
}
\end{gathered}
\end{equation}
and $\pi(\varphi) \in \overline{\mathbb M}^{(k)}$ belongs to the image of $(V^*/\rad V^*) \o (V/\rad V)$,
which is semisimple. Therefore the $\ggbar$-module $\overline{\mathbb M}_\lambda^{(k)}$ is semisimple,
and \eqref{eq:abstract semisimple quotients} implies that
$\mathbb M^{(k)} \subset \soc^k \mathbb M_\lambda$.

Next, let $\varphi \in \soc^k \mathbb M_\lambda$, and let $V$ be as in \eqref{eq:subspace generated by phi}. 
Then $ll(V) \le ll(\soc^k \mathbb M_\lambda) = k$ and $\varphi \in \mathbb M(V)$, hence $\varphi \in \mathbb M_\lambda^{(k)}$.
Therefore $\soc^k \mathbb M_\lambda \subset \mathbb M^{(k)}$, and the statement follows.
\end{proof}

\begin{prop}
\label{thm:M(lambda) composition series}
Let $\lambda \in -\mathbf P^{++}$. Then
\begin{equation}
\label{eq:M(P) character}
\ch \mathbb M_\lambda = \sum_{w \in W^\lambda} \ch \( M_{w\cdot\lambda}^* \o M_{w\cdot\lambda} \).
\end{equation}
Equivalently, for any $x,y \in W^\lambda$ we have
\begin{equation}
\label{eq:M(P) JH multiplicities}
[\mathbb M_\lambda : L_{x\cdot\lambda}^* \o L_{y\cdot\lambda}] = [P_{y\cdot\lambda}:L_{x\cdot\lambda}].
\end{equation}
\end{prop}

\begin{proof}
Let $x,y \in W^\lambda$, and let $L = L_{x\cdot\lambda}^* \o L_{y\cdot\lambda}$
be a simple constituent of the composition series of the $\ggbar$-module $\mathbb M_\lambda$, 
belonging to some layer $\overline{\mathbb M}_\lambda^{(k)}$ of the
filtration \eqref{eq:matrix elements filtration}.
It follows from \eqref{eq:diagram factors through} that $L$ 
can be generated by some $V \in \mathcal O_\lambda$ such that $ll(V) = k$.
We can assume without loss of generality that $V/\rad V \cong L_x$ and $V^*/ \rad V^* \cong L_y^*$.

Let $P = P_{\mu_1} \oplus \dots \oplus P_{\mu_j}$ be a projective cover of $V$. Then $V^* \subset \soc^k P^* $
and \eqref{eq:diagram factors through} implies that for $\mu_i \ne x\cdot\lambda$ the matrix elements of
$P_{\mu_i}$ do not contribute to $L_{x\cdot\lambda}^* \o L_{y\cdot\lambda} \subset \overline{\mathbb M}_\lambda^{(k)}$.
Therefore, $L$ can be generated by matrix elements of $P_{x\cdot\lambda}$, and moreover $L$
belongs to $\Phi(\soc^k P_{x\cdot\lambda}^* \o P_{x\cdot\lambda})$.

In other words, any simple constituent $L_{x\cdot\lambda}^* \o L_{y\cdot\lambda} \subset \overline{\mathbb M}_\lambda^{(k)}$
corresponds to a pair of submodules $L_{y\cdot\lambda} \subset \radl^1 P_{x\cdot\lambda}$ and 
$L_{x\cdot\lambda}^* \subset \socl^k P_{y\cdot\lambda}^* \cong (\radl^k P_{y\cdot\lambda})^*$. 
Since the projective module $P_{y\cdot\lambda}$ has simple top isomorphic to $L_{y\cdot\lambda}$, we see that
the occurrences of $L_{x\cdot\lambda}^* \o L_{y\cdot\lambda}$ in $\overline{\mathbb M}_\lambda^{(k)}$ are in bijective correspondence
with submodules $L_{x\cdot\lambda}^* \subset (\radl^k P_{y\cdot\lambda})^*$. Therefore,
\begin{equation}
[\mathbb M_\lambda : L_{x\cdot\lambda}^* \o L_{y\cdot\lambda}] = \sum_k
[\overline{\mathbb M}_\lambda^{(k)} : L_{x\cdot\lambda}^* \o L_{y\cdot\lambda}] = \sum_k
[(\radl^k P_{y\cdot\lambda})^* : L_{x\cdot\lambda}^*] = [P_{y\cdot\lambda} : L_{x\cdot\lambda}],
\end{equation}
which implies the desired statement.
\end{proof}

The simple constituents $L_{x\cdot\lambda}$ of the projective module $P_{y\cdot\lambda}$ are in bijective correspondence
with the morphisms from $\Hom_{\mathcal O}(P_{x\cdot\lambda},P_{y\cdot\lambda})$ - to any such morphism we associate
the constituent in $P_{y\cdot\lambda}$, determined by the image of the simple top of $P_{x\cdot\lambda}$. Therefore, we get

\begin{cor}
\label{thm:relation to homs}
For any $\lambda \in -\mathbf P^{++}$ and $x,y \in W^\lambda$ we have
\begin{equation}
[\mathbb M_\lambda: L_{x\cdot\lambda}^* \o L_{y\cdot\lambda}] = \dim \Hom_{\mathcal O}(P_{x\cdot\lambda},P_{y\cdot\lambda}).
\end{equation}
\end{cor}
Thus, the structure of the $\ggbar$-module $\mathbb M_\lambda$ is governed by the morphisms between projective modules
in the category $\mathcal O_\lambda$.

\begin{prop}
\label{thm:M(lambda) radical series}
For any $\lambda \in - \mathbf P^{++}$ the radical filtration of $\mathbb M_\lambda$ coincides with \eqref{eq:matrix elements filtration}.
\end{prop}
\begin{proof}
To verify that $\mathbb M_\lambda^{(k)} \subset \rad \mathbb M^{(k+1)}$, 
it suffices to show that for $k=1\dots,2l_\lambda$, every simple constituent of the layer $\overline{\mathbb M}_\lambda^{(k)}$ 
is non-trivially linked with the above layer $\overline{\mathbb M}_\lambda^{(k+1)}$.

Consider a simple component $L = L_{x\cdot\lambda}^* \o L_{y\cdot\lambda} \subset \overline{\mathbb M}_\lambda^{(k)}$, where $x,y \in W^\lambda$.
It follows from the proof of Proposition \ref{thm:M(lambda) composition series} that $L$ 
corresponds to a pair of simple constituents $L_{y\cdot\lambda} \in \radl^1 P_{y\cdot\lambda}$
and $L_{x\cdot\lambda}^* \in (\radl^{k} P_{y\cdot\lambda})^*$. Moreover, one can independently 
prove that $L$ can be represented by matrix elements of the big 
projective module $P_\lambda$, see Corollary \ref{thm:projectives suffice} below. Thus we may also assume
that $L$ is represented by a pair  $L_{y\cdot\lambda} \in \radl^j P_\lambda$
and $L_{x\cdot\lambda}^* \in (\radl^{j+k-1} P_\lambda)^*$. Since $k\le 2l_\lambda$, we either have
$j>1$ or $j+k-1<2l_\lambda+1$, or possibly both.

If $j>1$, then $\radl^{j-1} P_\lambda \ne 0$, and due to rigidity of $P_\lambda$ the constituent $L_{y\cdot\lambda}$
is non-trivially linked with some $L_{w\cdot\lambda} \subset \radl^{j-1} P_\lambda$. It follows that
$L_{x\cdot\lambda}^* \o L_{y\cdot\lambda}$ is non-trivially linked with $L_{x\cdot\lambda}^* \o L_{w\cdot\lambda} \subset \overline{\mathbb M}_\lambda^{(k+1)}$, and hence $L_{x\cdot\lambda}^* \o L_{y\cdot\lambda} \subset \rad \mathbb M_\lambda^{(k+1)}$.

If $j+k-1<2l_\lambda+1$, then $\radl^{j+k} P_\lambda \ne 0$, and as above we see that
$L_{x\cdot\lambda}^*$ is non-trivially linked with some $L_{w\cdot\lambda}^* \subset (\radl^{j+k} P_\lambda)^*$. Again, it follows that
$L_{x\cdot\lambda}^* \o L_{y\cdot\lambda}$ is non-trivially linked with $L_{w\cdot\lambda}^* \o L_{y\cdot\lambda} \subset \overline{\mathbb M}_\lambda^{(k+1)}$, and thus $L_{x\cdot\lambda}^* \o L_{y\cdot\lambda} \subset \rad \mathbb M_\lambda^{(k+1)}$.

We have established that for every $k$ we have $\mathbb M_\lambda^{(k)} \subset \rad \mathbb M^{(k+1)}$; the opposite
inclusion follows from \eqref{eq:abstract semisimple quotients}. Therefore, \eqref{eq:matrix elements filtration} 
is the radical filtration of $\mathbb M_\lambda$.
\end{proof}

Combining Proposition \ref{thm:M(lambda) socle series} and Proposition \ref{thm:M(lambda) radical series}, we get

\begin{cor}
For any $\lambda \in - \mathbf P^{++}$ the module $\mathbb M_\lambda$ is rigid and has Loewy length $2l_\lambda+1$.
\end{cor}

\section{Big projective modules and the structure of $\mathfrak R(G_0)$.}

\subsection{Matrix elements of big projectives modules}

\begin{prop}
\label{thm:over center}
For any $\lambda \in - \mathbf P^{++}$ there is a $\ggbar$-module isomorphism
\begin{equation}
\mathbb M(P_\lambda) \cong P_\lambda^* \o_{Z(\g)} P_\lambda
\end{equation}
\end{prop}

\begin{proof}
By construction, we have $\mathbb M(P_\lambda) = P_\lambda^* \o P_\lambda / \ker \Phi$. 
We need to check that $\ker \Phi$ coincides with the submodule $J \subset P_\lambda^* \o P_\lambda$, 
\begin{equation}
J = \biggr\<z^* v^* \o v  - v^* \o z \,  v \biggr\>, \qquad v \in P_\lambda, \, v^* \in P_\lambda^*, \, z \in Z(\g)
\end{equation}
It is easy to see that $J \subset \ker\Phi$; indeed, for every $v \in P_\lambda, \, v^* \in P_\lambda^*$ and $z \in Z(\g)$ we have
\begin{equation}
\Phi_{z^* v^* \o v} (x) = \<z^* v^*, x \, v\> = \<v^*, z \, x v\> =  \<v^*, x \, zv\> = \Phi_{v^* \o z v} (x), \qquad x \in \mathcal U(\g).
\end{equation}

To prove the opposite inclusion, we use Soergel's deformation of the projective modules.
We can think of it as a family of $\g$-modules $P_{\lambda;\eps}$, which are identical as vector spaces, but the action
of $\g$ depends on the deformation parameter $\eps \in \h^*$. 
The specialization $\eps = 0$ yields $P_{\lambda,\eps} \cong P_\lambda$; see \cite{So} for complete details.

When $\lambda$ is regular, we consider $\eps$ generic in a small neighborhood of $0 \in \h^*$. If $\lambda$
lies on some wall(s) of the Weyl chamber (i.e. $W_\lambda \ne \{e\}$), then we consider $\eps$ generic such that
$W_\eps = W_\lambda$, i.e. from the same wall(s) as $\lambda.$
For such $\eps$, the specialization $P_{\lambda;\eps}$ is a direct sum of the Verma modules,
\begin{equation}
\label{eq:deformed big projective}
P_{\lambda; \eps} \cong \bigoplus_{w\in W^\lambda} M_{w\cdot\lambda+\eps}.
\end{equation}
Moreover, the central characters corresponding to $w\cdot\lambda + \eps$ are all distinct, and therefore
\begin{equation}
\label{eq:deformed matrix elements}
P_{\lambda; \eps}^* \o_{Z(\g)} P_{\lambda;\eps} \cong 
\bigoplus_{w\in W^\lambda} M_{w\cdot\lambda+\eps}^* \o M_{w\cdot\lambda+\eps} \cong
\mathbb M(P_{\lambda; \eps}).
\end{equation}

Thus, generically we have $\ker \Phi_\eps = J_\eps$; the discrepancies between
$J_\eps$ and $\ker \Phi_\eps$ may occur when $\eps$ satisfies $\<\beta,\eps\> = 0$ for
one or more $\beta \in \Delta^+$. We are most interested, of course, in the extreme degenerate case $\eps = 0$.

As in \cite{So}, it suffices to check that the linear equations determining submodules $\ker \Phi_\eps$ and
$J_\eps$ still have the same rank in the {\it subgeneric} case, when $\eps$ satisfies $\<\beta,\eps\> = 0$ 
for a unique positive root $\beta$ such that $\<\beta,\lambda+\rho\> \ne 0$. This reduces our problem
to the verification of the Proposition for the Lie algebra $\mathfrak{sl}(2,\C)$, which we check directly
in the remaining part of the proof.

Let $\mu$ be a regular anti-dominant weight for $\mathfrak{sl}(2,\C)$; then $W^\mu = \{e,s\}$. 
Under the standard identification of the weight lattice for $\mathfrak{sl}(2,\C)$ with $\Z$, we have $\mu\in \{-2,-3,\dots\}$. 

We represent the structures of the big projective module $P_\mu$ and its dual $P_\mu^*$ pictorially, 
with the top blocks corresponding to the tops of the modules, and the bottom blocks representing the socles:

$$P_\mu \sim \begin{gathered}
\xymatrix@R=2pt{
*+[F-,]{L_\mu}\\ 
*+[F-,]{L_{s\cdot\mu}} \\
*+[F-,]{L_\mu}}
\end{gathered} \ , \qquad\qquad
P_\mu^* \sim \begin{gathered}
\xymatrix@R=2pt{
*+[F-,]{L_\mu^*}\\ 
*+[F-,]{L_{s\cdot\mu}^*} \\
*+[F-,]{L_\mu^*}
}\end{gathered} \ .
$$

Tensoring these two series, we get a filtration for $P_\mu^* \o P_\mu$, with layers depicted by
$$P_\mu^* \o P_\mu \sim \begin{gathered}
\xymatrix@R=2pt@C=-20pt{
&& *+[F-,]{L_\mu^* \o  L_\mu} &&\\ 
& *+[F-,]{L_{s\cdot\mu}^* \o  L_\mu} & \ar@{}|-\bigoplus & *+[F-,]{L_\mu^* \o L_{s\cdot\mu}} &\\
 *+[F-,]{L_\mu^* \o  L_\mu}  &  \ar@{}|-\bigoplus &  *+[F-,]{L_{s\cdot\mu}^* \o  L_{s\cdot\mu}}  &  \ar@{}|-\bigoplus &  *+[F-,]{L_\mu^* \o  L_\mu} \\
& *+[F-,]{L_\mu^* \o  L_{s\cdot\mu}} &  \ar@{}|-\bigoplus & *+[F-,]{L_{s\cdot\mu}^* \o L_\mu} &\\
&& *+[F-,]{L_\mu^* \o  L_\mu} &&
}\end{gathered}$$

The natural action of $Z(\g)$ on $P_\mu$ gives a surjection $Z(\g) \to \End_{\mathcal O}(P_\mu) \cong \C[Q]/\<Q^2\>$. 
In other words, the endomorphism algebra $\End_{\mathcal O}(P_\mu)$ is linearly spanned by the identity and an endomorphism $Q$, 
satisfying $Q^2 = 0$, which can be constructed as the composition map
$$Q: \quad \begin{gathered}
\xymatrix@R=2pt{
*+[F-,]{L_\mu} \ar@/^1pc/[rd] & & *+[F-,]{L_\mu} \\ 
*+[F-,]{L_{s\cdot\mu}} &  *+[F-,]{L_\mu} \ar@/_1pc/[rd] & *+[F-,]{L_{s\cdot\mu}}\\ 
*+[F-,]{L_\mu} & & *+[F-,]{L_\mu}}
\end{gathered} \  .$$
The action of $Z(\g)$ on $P_\mu^*$ is completely analogous.

It is now straightforward (cf. \cite{FS}) to see that both $J$ and $\ker \Phi$ are described by
$$J  = \ker \Phi \sim \begin{gathered}
\xymatrix@R=2pt@C=-5mm{
& *+[F-,]{L_\mu^* \o  L_\mu}  &  \\
*+[F-,]{L_\mu^* \o  L_{s\cdot\mu}} &  \ar@{}|-\bigoplus & *+[F-,]{L_{s\cdot\mu}^* \o L_\mu}\\
& *+[F-,]{L_\mu^* \o  L_\mu} &
}\end{gathered} \ ,$$
with the top constituent $L_\mu^* \o  L_\mu$ corresponding to the ``skew-symmetric'' part of the
middle layer of $P_\mu \o P_\mu^*$, which contained two copies of $L_\mu^* \o  L_\mu$.

Finally, there is only one singular antidominant weight $\mu = -1$, in which case the statement is
obvious because $P_{-1} = M_{-1} = L_{-1}$, and therefore $J = \ker \Phi = 0$.

This concludes the analysis of the $\mathfrak{sl}(2,\C)$ case, and the proof of the Proposition.
\end{proof}

In other words, we only need matrix elements of the big projective module $P_\lambda$ to span the entire block $\mathbb M_\lambda$.

\begin{cor}
\label{thm:projectives suffice}
For any $\lambda \in - \mathbf P^{++}$ we have $\mathbb M_\lambda = \mathbb M(P_\lambda)$.
\end{cor}

\begin{proof}
From the definitions it follows that $\mathbb M(P_\lambda) \subset \mathbb M_\lambda$. On the other hand, \eqref{eq:deformed matrix elements} implies,
$$\ch \mathbb M(P_\lambda) = \sum_{w\in W^\lambda} \ch \( M_{w\cdot\lambda}^* \o M_{w\cdot\lambda}  \),$$
or equivalently for any $x,y \in W^\lambda$ we have
\begin{multline}
[\mathbb M(P_\lambda):L_{x\cdot\lambda}^* \o L_{y\cdot\lambda}] = 
[\bigoplus_{w \in W^\lambda} M_{w\cdot\lambda}^* \o M_{w\cdot\lambda} : L_{x\cdot\lambda}^* \o L_{y\cdot\lambda}] = \\
= \sum_{w \in W^\lambda} [ M_{w\cdot\lambda}^* : L_{x\cdot\lambda}^*] [M_{w\cdot\lambda} : L_{y\cdot\lambda}] = 
\sum_{w \in W^\lambda} [ M_{w\cdot\lambda} : L_{x\cdot\lambda}] [P_{y\cdot\lambda} : M_{w\cdot\lambda}] =  [P_{y\cdot\lambda} : L_{x\cdot\lambda}].
\end{multline}
Comparing it with \eqref{eq:M(P) JH multiplicities}, we conclude that $\mathbb M(P_\lambda) = \mathbb M_\lambda$.
\end{proof}

\subsection{The $\ggbar$-module structure of $\mathfrak R(G_0)$ and the Peter-Weyl theorem}

Recall that the classical Peter-Weyl theorem asserts that for a compact Lie group $\bar G$ the matrix elements of finite-dimensional $\bar G$-modules 
produce an $L^2$-basis of the space $L^2(\bar G)$. Its algebraic version can be formulated
as follows:

\begin{thm}
Let $G$ be a simple complex Lie group, and let $\mathfrak R(G)$ denote the algebra of regular functions on $G$. 
Then there are isomorphisms of $G \times G$-modules
\begin{equation}
\mathfrak R(G) \cong \bigoplus_{\lambda \in \mathbf P^+} \mathbb M(L_\lambda) \cong  \bigoplus_{\lambda \in \mathbf P^+} L_\lambda \o L_\lambda^*
\end{equation}
\end{thm}

The first isomorphism $\mathfrak R(G) \cong \bigoplus_{\lambda \in \mathbf P^+} \mathbb M(L_\lambda)$ reflects the
spanning property, and the second isomorphism  $\mathfrak R(G) \cong \bigoplus_{\lambda \in \mathbf P^+} L_\lambda \o L_\lambda^*$
corresponds to the linear independence of matrix elements functions, corresponding to fixed bases of $L_\lambda$.

We can now formulate the projective analogue of the Peter-Weyl theorem.

\begin{thm}
\label{thm:projective Peter-Weyl}
Let $G$ be a simple complex group, and let $G_0$ be the big cell of $G$ associated with the Gauss decomposition.
Then there are $\ggbar$-module isomorphisms
\begin{equation}
\mathfrak R(G_0) \cong \bigoplus_{\lambda \in -\mathbf P^{++}} \mathbb M(P_\lambda) \cong 
\bigoplus_{\lambda \in -\mathbf P^{++}} P_\lambda \o_{Z(\g)} P_\lambda^*.
\end{equation}
\end{thm}

\begin{proof}
The map $\mathcal D: \mathcal U(\g) \to \Diff(G)$, induced by the identification of $\g$ with left-invariant vector
fields on $G$, yields a homomorphism of $\ggbar$-modules,
\begin{equation}
\vartheta: \mathfrak R(G_0) \to \mathcal U(\g)^*,\qquad
(\vartheta(\psi))(x) = (\mathcal D_x \psi)(e),
\end{equation}
where $e$ is the unit of the group $G$. We claim that $\vartheta$ is an injection.
Indeed, if $\psi \in \ker \vartheta$, then the function $\psi$ and all of its derivatives vanish at the unit of the group.
Since $\psi$ is regular, this means that $\psi$ is identically zero, and thus $\ker \vartheta = 0$.

It is obvious from the polynomial realization that $\mathfrak R(G_0)$ is locally $\n_- \o \n_+$-nilpotent,
and is $\h$-diagonalizable with respect to both regular actions. The same argument as in 
Proposition \ref{thm:matrix elements in O} implies that for any $\lambda \in - \mathbf P^{++}$ 
the submodule $\mathfrak R(G_0)_\lambda$ is isomorphic to a $\ggbar$-submodule of $\mathbb M_\lambda$.
On the other hand, combining \eqref{eq:R(G_0) JH multiplicities} and \eqref{eq:M(P) JH multiplicities}, we get
\begin{equation}
[\mathfrak R(G_0)_\lambda: L_{x\cdot\lambda} \o L_{y\cdot\lambda}^*] = 
[P_{x\cdot\lambda}: L_{y\cdot\lambda}] = [\mathbb M_\lambda: L_{x\cdot\lambda} \o L_{y\cdot\lambda}^*]
\end{equation}
for every $x,y \in W^\lambda$, which means that $\vartheta(\mathfrak R(G_0)_\lambda) = \mathbb M_\lambda$.
Taking the direct sum over $\lambda$ and using Propositions \ref{thm:over center} and \ref{thm:projectives suffice}, we get the desired statement.
\end{proof}

\subsection{Projective generators and their endomorphisms}
For any $\lambda \in -\mathbf P^{++}$, define
\begin{equation}
\mathcal P_\lambda = \bigoplus_{w\in W^\lambda} P_{w\cdot\lambda}
\end{equation}
to be the projective generator of $\mathcal O_\lambda$, and let $\mathcal A_\lambda$ denote its endomorphism algebra,
\begin{equation}
\mathcal A_\lambda = \End_{\mathcal O_\lambda} ( \mathcal P_\lambda ) = 
\bigoplus_{x,y\in W^\lambda}\Hom_{\mathcal O}(P_{x\cdot\lambda},P_{y\cdot\lambda}).
\end{equation}
The fundamental role of the algebra $\mathcal A_\lambda$ in the study of the category
$\mathcal O_\lambda$ is reflected by the equivalence of categories \cite{BGG}:
\begin{equation}
\mathcal O_\lambda \ \cong \ \text{finite dimensional } \mathcal A_\lambda\text{-mod}.
\end{equation}

\begin{thm}
\label{thm:Koszul Peter-Weyl}
For any $\lambda \in - \mathbf P^{++}$ there is an isomorphism of $\ggbar$-modules
\begin{equation}
\mathbb M_\lambda \cong \mathcal P_\lambda^* \o_{\mathcal A_\lambda} \mathcal P_\lambda.
\end{equation}
\end{thm}
\begin{proof}
It is clear that $\mathbb M_\lambda \cong (\mathcal P_\lambda^* \o \mathcal P_\lambda) / \ker \Phi$, 
because matrix elements of projectives in $\mathcal O_\lambda$ span $\mathbb M_\lambda$. 
We need to show that $\ker \Phi$ coincides with the $\ggbar$-submodule $\mathcal J_\lambda \subset \mathcal P_\lambda^* \o \mathcal P_\lambda$,
\begin{equation}
\label{eq:ideal by center}
\mathcal J_\lambda = \biggr\<\phi^*(v^*) \o v  - v^* \o \phi (v) \biggr\>, 
\qquad v \in \mathcal P_\lambda, \, v^* \in \mathcal P_\lambda^*, \, \phi \in \mathcal A_\lambda.
\end{equation}

Using the fact that $\phi \in \mathcal A_\lambda$ are $\g$-intertwining operators, we compute
\begin{equation}
\Phi_{\phi^*(v^*) \o v} (x) = \<\phi^*(v^*), x v\> = \<v^*, \phi (x v) \> =  
\<v^*, x \, \phi (v) \> = \Phi_{v^* \o \phi(v)} (x)
\end{equation}
for any $x \in \mathcal U(\g)$, which shows that $\mathcal J_\lambda \subset \ker \Phi$, and that $\mathbb M_\lambda$
is a quotient of $\mathcal P_\lambda^* \o_{\mathcal A_\lambda} \mathcal P_\lambda$. To 
complete the proof, we use Soergel's deformations of projective modules $P_{w\cdot\lambda}$ for $w\in W^\lambda$.

We use the same notion of generic $\eps$ as in the proof of Proposition \ref{thm:over center}.
Then for any $x \in W$ the specialization $P_{x\cdot\lambda;\eps}$ is a direct sum of the Verma modules,
\begin{equation}
P_{x\cdot\lambda; \eps} \cong \bigoplus_{y\in W^\lambda} M_{y\cdot\lambda+\eps} \o V_{x,y},
\end{equation}
where $V_{x,y}$ are the multiplicity spaces such that $\dim V_{x,y} = (P_{x\cdot\lambda}:M_{y\cdot\lambda}) = [M_{y\cdot\lambda}:L_{x\cdot\lambda}]$.
It follows that $\mathcal P_{\lambda;\eps}$ is isomorphic to the sum of Verma modules $M_{w\cdot\lambda+\eps}$ with positive multiplicities.
For generic $\eps$ the modules $M_{w\cdot\lambda+\eps}$ are irreducible, and the central characters corresponding to $w\cdot\lambda + \eps$ 
are all distinct. Therefore, elements of $\mathcal A_{\lambda;\eps} = \End(\mathcal P_{\lambda;\eps})$ just reshuffle the isotopic components 
$M_{w\cdot\lambda+\eps}$, and we have
\begin{equation}
\mathcal P_{\lambda; \eps}^* \o_{\mathcal A_{\lambda;\eps}} \mathcal P_{\lambda;\eps} \cong 
\bigoplus_{w\in W} M_{w\cdot\lambda+\eps}^* \o M_{w\cdot\lambda+\eps}.
\end{equation}

Thus, generically we have $\ker \Phi_\eps = \mathcal J_\eps$. As in the proof of Proposition \ref{thm:over center},
it suffices to check that no discrepancies occur when $\eps$ is subgeneric, which reduces the problem to the
rank one verification. The rest of the proof establishes the required statement for $\g = \mathfrak{sl}(2,\C)$.

As in the proof of Proposition \ref{thm:over center}, let $\mu \in \{-2,-3,\dots\}$ be a regular antidominant weight for
$\g$. The two indecomposable projective modules in $\mathcal O_\mu$ are the big projective module $P_\mu$ and the Verma
module $M_{s\cdot\mu}$. The algebra $\mathcal A_\mu$ has dimension 5. It has two obvious idempotents $\1_\mu$ and $\1_{s\cdot\mu}$,
and contains the inclusion $\iota: M_{s\cdot\mu} \to P_\mu$ and the map $\tau: P_\mu \to M_{s\cdot\mu}$, which factors through $L_\mu$;
we also set $\iota(P_\mu) = \tau(M_{s\cdot\mu}) = 0$.
The remaining element $Q \in \End_{\mathcal O}(P_\mu)$ is then equal to the composition $Q = \tau \circ \iota$. 

Using the idempotents $\1_\mu, \1_{s\cdot\mu}$, we immediately see that 
\begin{equation}
M_{s\cdot\mu}^* \o_{\mathcal A_\mu} P_\mu = P_\mu^* \o_{\mathcal A_\mu} M_{s\cdot\mu} = 0.
\end{equation}
Let now $v^* \o v \in M_{s\cdot\mu}^* \o M_{s\cdot\mu}$. We can find $u^* \in P_\mu^*$ such that $v^* = \iota^* (u^*)$, and therefore
\begin{equation}
v^* \o v = \iota^* (u^*) \o v  \equiv  u^* \o \iota(v) \in P_\mu^* \o P_\mu.
\end{equation}
Thus, any element in $\mathcal P_\mu^* \o_{\mathcal A_\mu} \mathcal P_\mu$ can be represented
by an element from $P_\mu^* \o P_\mu$, and Proposition \ref{thm:over center} implies
\begin{equation}
\mathcal P_\mu^* \o_{\mathcal A_\mu} \mathcal P_\mu \cong P_\mu^* \o_{\mathcal C_\mu} P_\mu \cong \mathbb M_\mu.
\end{equation}
This concludes the analysis of the $\mathfrak{sl}(2,\C)$ case, and the proof of the Proposition.
\end{proof}

The algebra $\mathcal A_\lambda$ is a graded algebra, and is the Koszul dual of the algebra $\Ext^\bullet(\mathcal L_\lambda,\mathcal L_\lambda)$, where
$\mathcal L_\lambda = \bigoplus_{w\in W^\lambda} L_{w\cdot\lambda}$, see \cite{So,BGS}. 
The algebra $\mathcal A_\lambda^\bullet$ is generated 
over $\mathcal A_\lambda^0$ by elements of $\mathcal A_\lambda^1$ with relations of degree 2. It is easy to see that
the degree 1 generators can be represented by elements of $\Hom_{\mathcal O}(P_{x\cdot\lambda},P_{y\cdot\lambda})$, such that
the simple top of $P_{x\cdot\lambda}$ is sent into $\radl^2 P_{y\cdot\lambda}$. More generally, elements of degree $k$ can 
be represented by morphisms sending the top of $P_{x\cdot\lambda}$ to the layer $\radl^{k+1} P_{y\cdot\lambda}$.
Thus, we can refine Corollary \ref{thm:relation to homs} as follows.

\begin{cor}
\label{thm:graded object}
Let $\lambda \in -\mathbf P^{++}$. Then the graded $\ggbar$-module $\operatorname{Gr}^\bullet \mathbb M_\lambda$, 
associated with the filtration \eqref{eq:matrix elements filtration}, is isomorphic to
\begin{equation}
\operatorname{Gr}^\bullet \mathbb M_\lambda \cong 
\bigoplus_{x,y \in W^\lambda} L_{x\cdot\lambda}^* \o  L_{y\cdot\lambda} \o \(\mathcal A_\lambda^\bullet\)_{x,y},
\end{equation}
where $(\mathcal A_\lambda^\bullet)_{x,y} = \Hom^\bullet_{\mathcal O_\lambda}( P_{x\cdot\lambda}, P_{y\cdot\lambda})$
are the multiplicity spaces, trivial as $\ggbar$-modules.
\end{cor}
Thus, $\mathbb M_\lambda$ is in a sense a ``categorification'' of the algebra $\mathcal A_\lambda$.
It would be interesting to interpret our results in the general context of mixed geometry \cite{BGS}.

\section{Generalizations to quantum groups.}

\subsection{Quantum groups with generic $q$}

The problem of constructing a $q$-deformation of the regular representation $\mathfrak R(G)$ was one
of the main motivations for the development of the theory of quantum groups. The quantum group 
$\mathcal U_q(\g)$ is a $q$-deformation of the Hopf algebra $\mathcal U(\g)$, generated by
$\mathbf E_i, \mathbf F_i, \mathbf K_i^{\pm1}$, subject to the standard relations, see e.g. \cite{L2}.
Using the analogy with the Lie algebraic realization of $\mathfrak R(G)$, we can define the quantum coordinate
algebra $\mathfrak R_q(G)$ as the Hopf dual $\mathcal U_q(\g)^*_{Hopf}$ of the quantum group, cf. \cite{L1,APW}.
The associative algebra $\mathfrak R_q(G)$ is equipped with two commuting quantum regular actions $\varrho_1, \varrho_2$.

When $q\in\C^\times$ is generic (or is regarded as a formal variable), the representation theory of $\mathcal U_q(\g)$ is parallel
to the classical case. In particular, one has the $q$-analogues of the BGG category $\mathcal O_q$, which contains the quantum 
counterparts $L_{\lambda,q}$, $M_{\lambda,q}$, $P_{\lambda,q}$
of the simple, Verma and projective modules.
The category $\mathcal O_q$ decomposes into direct sum of blocks $\mathcal O_{\lambda,q}$, according to the
characters of the center $Z_q(\g)$ of the quantum group.

We have the quantum version of the Peter-Weyl theorem, which asserts that
\begin{equation}
\mathfrak R_q(G) = \bigoplus_{\lambda \in \mathbf P^+} \mathbb M(L_{\lambda,q}) 
\cong \bigoplus_{\lambda \in \mathbf P^+} L_{\lambda,q}^* \o L_{\lambda,q},
\end{equation}
and it is also clear from this decomposition that the space of left $\mathcal U_q(\n_-)$-invariant elements
gives a model for finite-dimensional simple modules $L_{\lambda,q}$, yielding the quantum version of the Borel-Weil
realization.

The regular representation $\mathfrak R(G_0)$, studied in this paper, also has a quantum analogue.
Again, we define it by analogy with the Lie algebraic realization as a suitable subspace
of $\mathcal U_q(\g)^*$, so that we have $\mathfrak R_q(G_0) = \mathcal U_q(\n_-)^* \o \C[\mathbf P] \o \mathcal U_q(\n_+)^*$.
The space $\mathfrak R_q(G_0)$ is an associative, almost commutative algebra, and satisfies the
quantum version of the Peter-Weyl theorem for the big projective modules.

\begin{thm}
The algebra $\mathfrak R_q(G_0)$ has a decomposition
\begin{equation}
\mathfrak R_q(G_0) = \bigoplus_{\lambda \in -\mathbf P^{++}} \mathbb M(P_{\lambda,q}) 
\cong \bigoplus_{\lambda \in -\mathbf P^{++}} P_{\lambda,q}^* \o_{Z_q(\g)} P_{\lambda,q}.
\end{equation}
\end{thm}

The principal ingredient in the Borel-Weil realization of the big projective modules was 
the notion of the Whittaker vectors in $\mathfrak R(G_0)$. The obstacle to immediate
generalizations to the quantum case is the absence of nonsingular characters
$\boldsymbol \eta_q^+: \mathcal U_q(\n^+) \to \C$. 
Two possible approaches to quantum analogues of Whittaker functions were suggested in \cite{E,Se}, 
and were used in the study of the deformed quantum Toda system.

We choose the algebraic approach, based on the equivalence between Whittaker and Soergel
functors, which can also be generalized to the case when $q$ is a root of unity.
The analogue of Corollary \ref{thm:algebraic Whittaker} is given by

\begin{thm}
\label{thm:quantum algebraic Whittaker}
For any $\lambda \in - \mathbf P^{++}$, we have the isomorphisms
\begin{equation}
\Hom_{\varrho_1}(P_{\lambda,q}^*,\mathfrak R_q(G_0)) \cong P_{\lambda,q}, \qquad
\Hom_{\varrho_2}(P_{\lambda,q},\mathfrak R_q(G_0)) \cong P_{\lambda,q}^*.
\end{equation}
\end{thm}

The proofs of these quantum theorems are analogous to their Lie algebra counterparts.

\subsection{Quantum groups at roots of unity}

For simplicity, we will assume now that $\g$ is simply laced. Let $\ell\ge 3$ be an odd integer,
and let $q$ be a primitive $\ell$-th root of unity. Then $\mathcal U_q(\g)$ contains a central
ideal $\mathcal I$, generated by $E_i^\ell, F_i^\ell$ and $(K_i^\ell-1)$. The corresponding quotient
$\mathfrak U = \mathcal U_q(\g)/\mathcal I$ is a finite-dimensional Hopf algebra \cite{L1}.

We consider the non-semisimple category $\mathfrak O$ of finite-dimensional $\mathfrak U$-modules. 
The simple modules $L_{\lambda;\ell}$ and their
unique indecomposable projective covers $P_{\lambda;\ell}$ are parameterized by 
weights $\lambda$ from the restricted weight lattice $\mathbf P_\ell = \mathbf P/ (\ell \, \mathbf P)$, conveniently realized as a set by
\begin{equation}
\mathbf P_\ell \cong \left\{ \lambda \in \mathbf P^+ \, \bigr| \, \<\lambda,\a_i\> < \ell \text{ for all } \a_i \in \Pi \right\}.
\end{equation}

The Weyl group action on $\mathbf P$ induces an action of $W$ on $\mathbf P_\ell$, so that
\begin{equation}
w \circ \lambda =  w\cdot\lambda \mod \ell \, \mathbf P , \qquad  w \in W.
\end{equation}
For $\lambda \in \mathbf P_\ell$, we define its stabilizer $W_{\lambda;\ell}$ and the coset space $W^{\lambda;\ell}$
as in the Lie algebra case.

The linking principle implies that $\Ext^\bullet(L_{\lambda;\ell},L_{\mu;\ell}) = 0$ unless $\mu = w\circ\lambda$ for some $w \in W$.
Therefore, the category  $\mathfrak O$ decomposes into a direct sum of blocks, 
indexed by a set $\mathbf X_\ell$ of representatives of the Weyl group orbits in $\mathbf P_\ell$. 

The dual space $\mathfrak U^*$ is an associative algebra, and consists of the functionals from $\mathcal U_q(\g)^*$, 
vanishing on the ideal $\mathcal I$. It carries two commuting regular $\mathfrak U$-actions $\varrho_1$ and $\varrho_2$,
and as a $\mathfrak U \times \mathfrak U$-module decomposes into blocks
\begin{equation}
\mathfrak U^* \cong \bigoplus_{\lambda \in \mathbf X_\ell} \mathfrak U_\lambda^*.
\end{equation}

As in the Lie algebra case, the space $\mathfrak U^*_\lambda$ is spanned by the matrix elements of finite-dimensional $\mathfrak U^*$-modules.
However, a block of $\mathfrak U^*$ cannot be generated by a single indecomposable projective, and to give an analogue of the Peter-Weyl theorem
we must use the more general approach, provided by Theorem \ref{thm:Koszul Peter-Weyl}. We have

\begin{thm}
\label{thm:Peter-Weyl root of unity}
Let $\lambda \in \mathbf X_\ell$, and let $\mathcal P_{\lambda;\ell} = \bigoplus_{w \in W^{\lambda;\ell}} P_{w\circ\lambda;\ell}$ be the projective generator of
the block $\mathfrak U^*_\lambda$. Then $\mathfrak U^*_\lambda$ is spanned by matrix elements of $\mathcal P_{\lambda;\ell}$, and we have
\begin{equation}
\mathfrak U_\lambda^* = \mathcal P_{\lambda;\ell}^* \o_{\mathcal A_{\lambda;\ell}} \mathcal P_{\lambda;\ell}.
\end{equation}
\end{thm}

\begin{proof}
As in the proof of Proposition \ref{thm:matrix elements in O}, every functional every functional from $\mathfrak U^*$ 
can be interpreted as a matrix element of a finite-dimensional $\mathfrak U$-module, and since $\mathfrak O$ has enough
projectives (see \cite{APW}), it suffices to consider matrix elements of $\mathcal P_{\lambda;\ell}$.

For $\lambda\in\mathbf P_\ell$ the projective $\mathfrak U$-module, remains projective when regarded as a $\mathcal U(\g)$-module.
We consider the $q$-version of Soergel's deformation $P_{\lambda+\eps;\ell}$ of the projective $\mathcal U_q(\g)$-modules in $\mathfrak O$, 
and corresponding matrix elements. Generically, the deformed projectives split into a direct sum of simple modules, and
the statement is clear. To check that the spaces still coincide when $\lambda$ becomes integral, it suffices again to consider subgeneric 
$\eps$, which reduces the problem to the rank one case.

Thus, we assume that $\g= \mathfrak {sl}(2,\C)$. The regular orbits of the affine Weyl group are $\{\mu,\mu'\}$, where
$\mu \in \mathbf X_\ell = \{0,1,\dots,\frac{\ell-3}2\}$ and $\mu' = \ell - \mu - 2$.
The corresponding projective modules are depicted by
$$P_{\mu;\ell} \ \sim \ \begin{gathered}
\xymatrix@R=2pt@C=-3pt{
& *+[F-,]{L_{\mu;\ell}}  &\\
*+[F-,]{L_{\mu';\ell}} \ar@{}[rr]|-\bigoplus && *+[F-,]{L_{\mu';\ell}} \\
& *+[F-,]{L_{\mu;\ell}} & 
}\end{gathered} \  , \qquad
P_{\mu';\ell} \ \sim \ \begin{gathered}
\xymatrix@R=2pt@C=-5pt{
& *+[F-,]{L_{\mu';\ell}}  &\\
*+[F-,]{L_{\mu;\ell}} \ar@{}[rr]|-\bigoplus && *+[F-,]{L_{\mu;\ell}} \\
& *+[F-,]{L_{\mu';\ell}} & 
}\end{gathered} \  .
$$
Thus, $\mathcal P_{\mu;\ell} = P_{\mu;\ell} \oplus P_{\mu';\ell}$, and it is easy to
describe the graded algebra $\mathcal A_{\mu;\ell} = \End_{\mathfrak U}(\mathcal P_{\mu;\ell})$ of dimension eight.
Indeed, the space $\End(P_{\mu;\ell})$ contains the identity map of degree 0, and the nontrivial
endomorphism of degree 2, factoring through $L_\mu$ . The space $\Hom(P_{\mu;\ell},P_{\mu';\ell})$
has two generators of degree 1, arising from the ``baby Verma'' flags of the projective modules.
Reversing the role of $\mu$ and $\mu'$, we get the remaining four generators of $\mathcal A_{\mu;\ell}$.

The verification of the fact that the analogue of \eqref{eq:ideal by center} coincides with the
kernel of the matrix elements map is a straightforward exercise, which we leave to the reader.
\end{proof}

The obvious modification of Corollary \ref{thm:graded object} describes the structure of 
the regular block $\mathfrak U^*_\mu$. For $\g=\mathfrak{sl}(2,\C)$ it is given by
\begin{equation}
\label{eq:Anechka diagram}
\mathfrak U_\mu^* \quad\sim\quad \begin{gathered}
\xymatrix@R=2pt@C=-4pt{
&& *+[F-,]{L_{\mu';\ell}^* \o  L_{\mu';\ell}}  \ar@{}[rr]|-\bigoplus & \quad & *+[F-,]{L_{\mu;\ell}^* \o L_{\mu;\ell}} &\\
 *+[F-,]{L_{\mu';\ell}^* \o  L_{\mu;\ell}}  \ar@{}[rr]|-\bigoplus & \quad & *+[F-,]{L_{\mu;\ell}^* \o  L_{\mu';\ell}}  \ar@{}[rr]|-\bigoplus & \quad  &   
 *+[F-,]{L_{\mu';\ell}^* \o  L_{\mu;\ell}}  \ar@{}[rr]|-\bigoplus & \quad
&  *+[F-,]{L_{\mu;\ell}^* \o  L_{\mu';\ell}}\\
&& *+[F-,]{L_{\mu;\ell}^* \o L_{\mu;\ell}}  \ar@{}[rr]|-\bigoplus & \quad & *+[F-,]{L_{\mu';\ell}^* \o  L_{\mu';\ell}} &
}\end{gathered} \quad,
\end{equation}
and was originally obtained in \cite{AGL} by explicit computations. A representation-theoretic derivation of \eqref{eq:Anechka diagram}
using matrix elements of projective modules was done by A. Lyakhovskaya (A. Lachowska) and the author.

The analogue of Corollary \ref{thm:algebraic Whittaker} and Theorem \ref{thm:quantum algebraic Whittaker} holds for $\mathfrak U^*$, 
and can also be proved using the self-duality of $\mathfrak U^*$.



\begin{thebibliography}{}
\frenchspacing

\bibitem[AGL]{AGL}
A. Alekseev, D. Glushchenkov, A. Lyakhovskaya,
Regular representation of the quantum group ${\mathfrak sl}_q(2)$ ($q$ is a root of unity).
Algebra i Analiz 6 (1994), no. 5, 88--125.

\bibitem[APW]{APW}
H. Andersen, P. Polo, K. Wen,
Representations of quantum algebras.  Invent. Math.  104  (1991),  no. 1, 1--59. 

\bibitem[Ba]{Ba}
E. Backelin,
Representation of the category $\mathcal O$ in Whittaker categories.  Internat. Math. Res. Notices  1997,  no. 4, 153--172. 

\bibitem[Be]{Be}
J. Bernstein,
Trace in categories.
Operator algebras, unitary representations, enveloping algebras, and invariant theory (Paris, 1989), 417--423,
Progr. Math., 92,
Birkhäuser Boston, 1990. 

\bibitem[BGS]{BGS}
A. Beilinson, V. Ginzburg, W. Soergel,
Koszul duality patterns in representation theory. 
J. Amer. Math. Soc. 9 (1996), no. 2, 473--527.

\bibitem[BGG]{BGG}
J. Bernstein, I. Gelfand, S. Gelfand,
A certain category of $\g$-modules.
(Russian) Funkcional. Anal. i Prilo\v zen. 10 (1976), no. 2, 1--8. 


\bibitem[Bo]{Bo}
A. Borel,
Linear representations of semi-simple algebraic groups.
Proc. Sympos. Pure Math., Vol. 29, Amer. Math. Soc. (1975), 421--439. 

\bibitem[Di]{Di}
J. Dixmier, Alg\`ebres enveloppantes. 
Gauthier-Villars \'Editeur, Paris, 1974.


\bibitem[E]{E}
P. Etingof,
Whittaker functions on quantum groups and $q$-deformed Toda operators.
Differential topology, infinite-dimensional Lie algebras, and applications, 9--25,
Amer. Math. Soc. Transl. Ser. 2, 194, 1999. 


\bibitem[FS]{FS}
I. Frenkel, K. Styrkas,
Modified regular representations of affine and Virasoro algebras, VOA structure and semi-infinite cohomology.
math.QA/0409117.


\bibitem[Ir]{Ir}
R. Irving,
Projective modules in the category ${\mathcal O}\sb S$: Loewy series.  Trans. Amer. Math. Soc.  291  (1985),  no. 2, 733--754. 


\bibitem[K1]{K}
B. Kostant,
On Whittaker vectors and representation theory.
Invent. Math. 48 (1978), no. 2, 101--184.

\bibitem[K2]{KToda}
B. Kostant,
The solution to a generalized Toda lattice and representation theory.  Adv. in Math.  34  (1979), no. 3, 195--338.

\bibitem[L1]{L1}
G. Lusztig,
Finite-dimensional Hopf algebras arising from quantized universal enveloping algebra.
J. Amer. Math. Soc.  3  (1990),  no. 1, 257--296. 

\bibitem[L2]{L2}
G. Lusztig,
Introduction to quantum groups. 
Progress in Mathematics, 110. Birkhäuser Boston, 1993.

\bibitem[Se]{Se}
A. Sevostyanov,
Quantum deformation of Whittaker modules and the Toda lattice.
Duke Math. J. 105 (2000), no. 2, 211--238.


\bibitem[So]{So}
W. Soergel,
Kategorie $\mathcal O$, perverse Garben und Moduln Žüber den Koinvarianten zur Weylgruppe. J. Amer. Math. Soc.  3  (1990),  no. 2, 421--445.

\bibitem[WW]{WW}
E..T. Whittaker, G.N. Watson,
A course in modern analysis.
London, Cambridge University Press, 1963.


\end{thebibliography}
\end{document}